\documentclass[11pt]{article}

\usepackage{pifont}
\usepackage{mathrsfs}
\usepackage{amscd}
\usepackage{amsmath}
\usepackage{latexsym}
\usepackage{amsfonts}
\usepackage{amssymb}
\usepackage{amsthm}
\usepackage{graphicx}
\usepackage{amsmath,amsfonts,ams}
 \newtheorem{proposition}{Proposition}[section]
 \newtheorem{definition}{Definition}[section]
 \newtheorem{lemma}{Lemma}[section]
 \newtheorem{theorem}{Theorem}[section]
 \newtheorem{corollary}{Corollary}[section]
 \newtheorem{remark}{Remark}[section]

\makeatletter
\newcommand{\Extend}[5]{\ext@arrow0099{\arrowfill@#1#2#3}{#4}{#5}}
\makeatother

\begin{document}

 \title{ \bf  Energy Scattering for a  Klein-Gordon  Equation with a Cubic Convolution}
 \author{{Changxing Miao$^1$, \  Jiqiang Zheng$^2$  }\\
        {\small  $^1$ Institute of Applied Physics and Computational Mathematics}\\
        {\small P. O. Box 8009,\ Beijing,\ China,\ 100088;}\\
      {\small $^2$ The Graduate School of China Academy of Engineering Physics  }\\
   {\small P. O. Box 2101,\ Beijing,\ China,\ 100088 ;} \\
   {\small (miao\_changxing@iapcm.ac.cn; \; zhengjiqiang@gmail.com ) }\\
    \date{ }
       }

\maketitle

\begin{abstract}
In this paper, we study the global well-posedness and scattering
problem in the energy space for both focusing and defocusing the Klein-Gordon-Hartree equation in
the spatial dimension $d \geq 3$.  The main difficulties are the
absence of an interaction Morawetz-type estimate and of a Lorentz
invariance which enable one to control the momentum. To compensate, we utilize the strategy derived from concentration compactness ideas, which was first introduced by Kenig and Merle \cite{KM} to the scattering problem.
Furthermore, employing technique from \cite{Pa2}, we consider a virial-type
identity in the direction orthogonal to the momentum vector so as to
control the momentum in the defocusing case. While in the focusing case, we
show that the scattering holds when the initial data $(u_0,u_1)$ is
radial, and the energy $E(u_0,u_1)<E(\mathcal{W},0)$ and $\|\nabla
u_0\|_2^2+\|u_0\|_2^2<\|\nabla \mathcal{W}\|_2^2+\|\mathcal{W}\|_2^2$, where $\mathcal{W}$ is the ground state.
\end{abstract}

\begin{center}
 \begin{minipage}{120mm}
 \end{minipage}
\end{center}


\section{Introduction}
\setcounter{section}{1}\setcounter{equation}{0} This paper is
devoted to the study of the Cauchy problem of the Klein-Gordon-Hartree equation
\begin{align} \label{equ1}
\begin{cases}    \ddot{u} - \Delta u  +  u +  f(u)=  0,  \qquad  (t,x) \in
\mathbb{R}\times\mathbb{R}^d, d \geq 3,\\
u(0,x)=u_0(x),~u_t(0,x)=u_1,\end{cases}
\end{align}
where $f(u)=\mu(|x|^{-\gamma}*|u|^2) u $, $2<\gamma<\min(4,d)$,
$\mu=\pm 1$ with $\mu=1$ known as the defocusing case and $\mu=-1$
as the focusing case. Here $u$ is a real-valued function defined in
$\mathbb{R}^{d+1}$, the dot denotes the time derivative, $\Delta$ is
the Laplacian in $\mathbb{R}^{d}$, $|x|^{-\gamma}$ is called the
potential, and $*$ denotes the spatial convolution in
$\mathbb{R}^{d}$.

Formally, the solution $u$ of \eqref{equ1} conserves  the energy
\begin{equation*}\label{econ}
\aligned E(u(t),\dot{u}(t))=&\frac12 \int_{\mathbb{R}^d}
\big(\big|\dot{u}(t,x) \big|^2 +\big| \nabla u(t,x) \big|^2 + \big| u(t,x) \big|^2
\big)dx\\& + \frac{\mu}{4} \iint_{\mathbb{R}^d\times\mathbb{R}^d}
 \frac{|
u(t,x)|^2 |u(t,y)|^2}{|x-y|^\gamma} dxdy\\
\equiv&E(u_0,u_1),
\endaligned
\end{equation*}
and the momentum
\begin{equation}
P(u)(t)=\int_{\mathbb{R}^d}u_t(t,x)\nabla u(t,x) dx=P(u)(0).
\end{equation}

The scattering theory for the Klein-Gordon equation with
$f(u)=\mu|u|^{p-1}u$  has been intensively studied in \cite{Br84},
\cite{Br85},  \cite{GiV85b}, \cite{IMN}, \cite{Na99b} and
\cite{Na01}. For $\mu=1$ and
\begin{equation}\label{}
   1+\frac{4}{d}<p<1+\gamma_d\frac{4}{d-2},
\quad
 \gamma_d=\left\{ \aligned
    &1 ,& 3\leq d\leq 9; \\
     &\frac{d}{d+1},& d\geq 10.
\endaligned
\right.
\end{equation}
Brenner \cite{Br85} established the scattering results in the energy
space, which does not contain all subcritical cases for $d\geq 10$.
Thereafter, Ginibre and Velo \cite{GiV85b} exploited the
Birman-Solomjak space $\ell^m(L^q,I,B)$ in \cite{BiS75} and the
delicate estimates to improve the results in \cite{Br85}, which
covered all subcritical cases. Finally K. Nakanishi \cite{Na99b}
obtained the scattering results for the critical case by the
strategy of induction on energy  \cite{CKSTT07} and a new
Morawetz-type estimate. And recently, S. Ibrahim, N. Masmoudi and
K. Nakanishi\cite{IMN, IMN1} utilized the concentration compactness
ideas to give the scattering threshold for the focusing nonlinear
Klein-Gordon equation. Their method also works for the defocusing case.

On the other hand, the scattering theory for the Hartree equation
$$i\dot{u}=-\Delta u+\mu(|x|^{-\gamma}*|u|^2)u$$
has been also  studied by many authors (see
\cite{GiV00,LiMZ09,MXZ07a,MXZ09,MXZ09a,MXZ09b,MXZ10,MXZ07b}). For
the subcritical defocusing case, Ginibre and Velo
\cite{GiV00} derived the associated Morawetz inequality and
extracted an useful Birman-Solomjak type estimate to obtain the
asymptotic completeness in the energy space. Nakanishi \cite{Na99d}
improved the results by a new Morawetz estimate which doesn't depend
on the nonlinearity. For the critical case, Miao, Xu and Zhao
\cite{MXZ07a} took advantage of a new kind of the localized Morawetz
estimate, which is also independent of the nonlinearity, to rule out
the possibility of the energy concentration at origin and
established the scattering results in the energy space for the
radial data in dimension $d\geq 5$. We refer also to
\cite{MXZ09,MXZ09a,MXZ09b,MXZ10,MXZ07b} for the general data and
focusing case.

For the equation $(\ref{equ1})$, using the ideas of Strauss
\cite{St81a}, \cite{St81b}, Pecher \cite{Pe85} and Mochizuki
\cite{Mo89} showed that if $d \geq 3$, $2\leq \gamma < \min(d, 4)$,
then global well-posedness and scattering results with small
data hold in the energy space $H^1\times L^2$. We refer also  to
Miao-Zhang\cite{MZ} where the low regularity for the cubic
convolution defocusing Klein-Gordon-Hartree equation is discussed.
In this paper, we develop in the energy space a complete scattering theory for
$(\ref{equ1})$  with the subcritical nonlinearity under
some suitable assumptions.  Compared with the classical Klein-Gordon
equation with the local nonlinearity $f=|u|^{p-1}u$, the
nonlinearity $f(u)=\mu(|x|^{-\gamma}*|u|^2) u $ is nonlocal, which brings us many
difficulties. The main difficulty is the absence of a Lorentz
invariance which could be used to control the momentum
efficiently. We will overcome this difficulty by considering a
Virial-type identity in the direction orthogonal to the momentum
vector following the technique in
\cite{Pa2} for the defocusing case. Unfortunately, the type of
virial identify in \cite{Pa2} can not work for the focusing case.
Inspired by the method in \cite{IMN} and \cite{MW}, we get over this
difficulty by some new variational framework and the profile
decomposition under the restriction that the initial data is radial.

Now we state our first result

\begin{theorem}\label{theorem}
Assume that $\mu=1$, $d\geq 3$, $2<\gamma<\min\{d,4\}$ and
$(u_0,u_1)\in H^1(\mathbb{R}^d)\times L^2(\mathbb{R}^d)$. Then there
exists a unique global solution $u(t)$ of \eqref{equ1} which
scatters in the sense that
 there exist  solutions $v_\pm$ of the free Klein-Gordon equation
\begin{equation}\label{le}
    \ddot{v} - \Delta v  + v =  0
\end{equation}
with $(v_\pm(0), \dot{v}_\pm(0))\in H^1\times L^2$ such that
\begin{equation}\label{1.2}
\big\|\big(u(t), \dot{u}(t)\big)-\big(v_\pm(t),
\dot{v}_\pm(t)\big)\big\|_{H^1\times L^2} \longrightarrow 0,\quad
\text{as}\quad t\longrightarrow \pm\infty.\end{equation}
\end{theorem}

In the focusing case, one can not expect to establish a similar
result as Theorem \ref{theorem} without any other restriction. In
fact, let $\mathcal{W}$ be an element of the ground states which satisfy the
elliptic equation
$$-\Delta\psi+\psi=(|x|^{-\gamma}*|\psi|^2)\psi.$$
One may find that $\mathcal{W}$ is a non-scattering solution of \eqref{equ1}
with the finite energy $E(\mathcal{W},0)$. The existence of ground state $\mathcal{W}$
was proved by \cite{Liu}. We will discuss the ground state in the
Section 3.

Now we state our second result as follows

\begin{theorem}\label{theorem1.2}
Assume $d\geq3,\ \mu=-1,$ and $2<\gamma<\min\{d,4\}$. Let
$(u_0,u_1)\in H^1\times L^2$ be radial with
$$E(u_0,u_1)<E(\mathcal{W},0),$$
and u be the corresponding solution of \eqref{equ1} with maximal
interval of existence $I=(-T_-(u_0,u_1),T_+(u_0,u_1))$.

(i) If $\|\nabla u_0\|_2^2+\|u_0\|_2^2<\|\nabla \mathcal{W}\|_2^2+\|\mathcal{W}\|_2^2$,
then $u$ is global and scatters.

(ii) If $\|\nabla u_0\|_2^2+\|u_0\|_2^2>\|\nabla \mathcal{W}\|_2^2+\|\mathcal{W}\|_2^2$,
then $u$ blows up both forward and backward in finite time, i.e.
$T_-(u_0,u_1),T_+(u_0,u_1)<\infty.$
\end{theorem}

\begin{remark}
The case $\|\nabla u_0\|_2^2+\|u_0\|_2^2=\|\nabla \mathcal{W}\|_2^2+\|\mathcal{W}\|_2^2$
is not compatible with $E(u_0,u_1)<E(\mathcal{W},0)$ by Lemma \ref{mini1} and
Lemma \ref{equi} below. In this sense, we give a complete
classification for solutions under the restriction
$E(u_0,u_1)<E(\mathcal{W},0)$.
\end{remark}

The outline for the proof of Theorem \ref{theorem} and Theorem
\ref{theorem1.2}: first, we define the scattering size of a solution
to \eqref{equ1} on a time interval $I\subset\mathbb{R}$ by
\begin{equation}\label{ST}
ST(I)=[W](I)\cap[K](I),\end{equation} where
$$[W](I)=L_t^{\frac{2(d+1)}{d-1}}(I;B^{\frac12}_{\frac{2(d+1)}{d-1},2}(\mathbb{R}^d)),\
[K](I)=L^{\frac{2(d+2)}{d}}(I;B^{\frac12}_{\frac{2(d+2)}{d},2}(\mathbb{R}^d)).$$
Then it is easy to see (cf Proposition \ref{prop2.3} below) that,
 scattering for the Klein-Gordon-Hartree
equation is implied by the finiteness of the scattering size $\Vert
u\Vert_{ST(I)}$ with $I=\mathbb{R}$.

Second, we define the function $\Lambda$ by
\begin{equation}\label{znorm}
\Lambda(E)=\sup\{\Vert u\Vert_{ST(\mathbb{R})}:E(u,u_t)\le E\}
\end{equation}
where the supremum is taken over all nonlinear solutions of
\eqref{equ1} with energy not greater than $E$, and define
\begin{align*}\label{emax}E_{max}^+&=\sup\{E:\Lambda(E)<+\infty\},\\
E_{max}^-&=\sup\{E:\Lambda(E)<+\infty, \|\nabla
u_0\|_2^2+\|u_0\|_2^2<\|\nabla \mathcal{W}\|_2^2+\|\mathcal{W}\|_2^2\},
\end{align*}
and \begin{align*}E_{max}=\begin{cases} E_{max}^+,\ \text{if}\
\mu=1,\\
E_{max}^-,\ \text{if}\ \mu=-1.\end{cases}\end{align*} In particular,
by Lemma \ref{equi},  $E_{max}^-$ is equivalent to
$$E_{max}^-=\sup\{E:\Lambda(E)<+\infty, K_1(u_0)>0\}.$$
Our goal next is to prove that $E_{max}^+=+\infty$ and
$E_{max}^-=E(\mathcal{W},0)$. We argue by contradiction. We show that if
$E_{max}^+<+\infty$ (or $E_{max}^-<E(\mathcal{W},0)$), then there exists a
nonlinear solution of \eqref{equ1} with energy be exactly $E_{max}$. Moreover, this solution satisfies some strong compactness properties. This is
completed in Section 6 where we utilize the profile decomposition
that was established in \cite{IMN}, and a strategy introduced by
Kenig and Merle \cite{KM}. We consider a virial-type identity in the
direction orthogonal to the momentum vector  following the technique \cite{Pa2} in the defocusing
case to obtain a contradiction. We refer to Section 7 for more details.

The paper is organized as follows. In Section $2$, we deal with the local theory
 for the equation $(\ref{equ1})$. In Section 3, we
discuss the property of the ground state. In Section 4, we prove the
blow up part of Theorem \ref{theorem1.2}. In Section $5$, we give
the linear and nonlinear profile decomposition and show some
properties of the profile. In Section 6, we extract a critical
solution. Finally in Section $7$, we preclude the critical solution,
which completes the proof of Theorem \ref{theorem} and Theorem
\ref{theorem1.2}.

We conclude the introduction by giving some notations which
will be used throughout this paper. We always assume the spatial
dimension $d\geq 3$ and let $2^*=\frac{2d}{d-2}$. For any $r, 1\leq
r \leq \infty$, we denote by $\|\cdot \|_{r}$ the norm in
$L^{r}=L^{r}(\mathbb{R}^d)$ and by $r'$ the conjugate exponent
defined by $\frac{1}{r} + \frac{1}{r'}=1$. For any $s\in
\mathbb{R}$, we denote by $H^s(\mathbb{R}^d)$ the usual Sobolev
space. Let $\psi\in \mathcal{S}(\mathbb{R}^d)$ be such that
$\text{supp}\ {\widehat{\psi}} \subseteq \big\{\xi: \frac{1}{2}
\leq|\xi| \leq 2 \big\}$ and $ \sum_{j\in \mathbb{Z}} \widehat{\psi}
(2^{-j} \xi) = 1 $ for $\xi \neq 0.$ Define $\psi_0$ by
$\widehat{\psi}_0 = 1 -
 \sum_{j\geq 1} \widehat{\psi} (2^{-j} \xi).$  Thus $\text{supp}\
\widehat{\psi}_0 \subseteq \big\{\xi: |\xi| \leq 2 \big\}$ and
$\widehat{\psi}_0 = 1$ for $|\xi| \leq 1$. We denote by $\Delta_j$
and $\mathcal{P}_0$ the convolution operators whose symbols are
respectively given by $\widehat{\psi}(\xi/2^{j})$ and
$\widehat{\psi}_0(\xi)$. For $s \in \mathbb{R}, 1\leq r \leq
\infty$, the inhomogeneous Besov space $ B^{s}_{r, 2}(\mathbb{R}^d)$
is defined by
$$ B^{s}_{r, 2}(\mathbb{R}^d) = \bigg\{ u \in
\mathcal{S}'(\mathbb{R}^d), \|\mathcal{P}_0 u\|^2_{L^r}+
\big\|2^{js} \|\Delta_j u\|_{L^r} \big\|^2_{l^2_{j\in \mathbb{N}}} <
\infty \bigg\}.$$ For details of Besov space, we refer to \cite{BL76}. For any interval
$I\subset \mathbb{R}$ and any Banach space $X$ we denote by ${\cal
C}(I; X)$ the space of strongly continuous functions from $I$ to $X$
and by $L^q(I; X)$ the space of strongly measurable functions from
$I$ to $X$ with $\|u(\cdot); X\|\in L^q(I).$ Given $d,$ we define,
for $2\le r\le \infty$,
$$\delta(r)=d\Big(\frac12-\frac1{r}\Big).$$
Sometimes abbreviate $\delta(r)$, $\delta(r_i) $ to
$\delta,~\delta_i $ respectively. We denote by $\langle\cdot,
\cdot\rangle$ the scalar product in $L^2$. We let $L_*^p$ denote the
weak $L^p$ space.

\section{Preliminaries}
 \setcounter{section}{2}\setcounter{equation}{0}

 \subsection{Strichartz estimate}
In this section, we consider the Cauchy problem  for the equation
$(\ref{equ1})$
\begin{equation} \label{equ2}
    \left\{ \aligned &\ddot{u} - \Delta u  +  u +  f(u)=  0, \\
    &u(0)=u_0,~\dot{u}(0)=u_1.
    \endaligned
    \right.
\end{equation}
The integral equation for the Cauchy problem $(\ref{equ2})$ can be
written as
\begin{equation}\label{inte1}
u(t)=\dot{K}(t)u_0 + K(t)u_1-\int^{t}_{0}K(t-s)f(u(s))ds,
\end{equation}
or
\begin{equation}\label{inte2}
{u(t)\choose \dot{u}(t)} = V_0(t){u_0(x) \choose u_1(x)}
-\int^{t}_{0}V_0(t-s){0 \choose f(u(s))} ds,
\end{equation}
where
$$K(t)=\frac{\sin(t\omega)}{\omega}, \quad V_0(t) = {\dot{K}(t), K(t)
\choose \ddot{K}(t), \dot{K}(t)}, \quad \omega=\big( 1-\Delta\big)^{1/2}.$$

Let $U(t)=e^{it\omega}$, then
\begin{equation*}
\dot{K}(t)= \frac{U(t)+U(-t)}{2}, \qquad  K(t)=
\frac{U(t)-U(-t)}{2i\omega}.
\end{equation*}

We begin by recalling the definition of  strong solution to the
Cauchy problem.
\begin{definition}\label{strong}
We call $u$ a strong solution to \eqref{equ2} on a time interval
$I$, if $u\in C(I;H^1)\cap C^1(I;L^2)$ and satisfies the Duhamel
formula \eqref{inte1} in the sense of tempered distributions for
every $t\in I.$
\end{definition}
Now we recall the following dispersive estimate for the operator
$U(t)=e^{it\omega}$.
\begin{lemma}[\cite{Br85}, \cite{GiV85b}]\label{lem21}
Let $2\leq r\leq \infty$ and $0\leq \theta\leq 1$. Then
\begin{equation*}
\big\|e^{i\omega t}f
\big\|_{B^{-(d+1+\theta)(\frac12-\frac1r)/2}_{r, 2}} \leq \mu(t)
\big\|f\big\|_{B^{(d+1+\theta)(\frac12-\frac1r)/2}_{r', 2}},
\end{equation*}
where
\begin{equation*}
\mu(t)=C \min\bigg\{ |t|^{-(d-1-\theta)(\frac12-\frac1r)_{+} },
|t|^{-(d-1+\theta)(\frac12-\frac1r)}\bigg\}.
\end{equation*}
\end{lemma}

According to the above lemma, the abstract duality and interpolation
argument(see \cite{GiV95}, \cite{KeT98}), we have the following
Strichartz estimates.
\begin{lemma}[\cite{Br85},\cite{GiV85b},\cite{MZF}]\label{lem22}
Let $0\leq\theta_i \leq 1$, $\rho_i \in \mathbb{R}$, $2\leq q_i, r_i
\leq \infty,~i=1,2$. Assume that $(\theta_i,d,q_i,r_i)\neq (0,3,2,\infty)$
satisfy the following admissible conditions
\begin{equation}
\left\{ \aligned \label{rl} 0\leq \frac{2}{q_i} &\leq
\min\Big\{(d-1+\theta_i)(\frac{1}{2}-\frac{1}{r_i}),
1\Big\},~~~i=1,2
   \\
&\rho_1+(d+\theta_1)(\frac{1}{2}-\frac{1}{r_1})-\frac{1}{q_1} =\mu,
\\&\rho_2+(d+\theta_2)(\frac{1}{2}-\frac{1}{r_2})-\frac{1}{q_2}
=1-\mu.
\endaligned\right.
\end{equation}
Then, for $f \in H^\mu$, we have
\begin{align}\label{str1}
\big\| U(\cdot) f\big\|_{L^{q_1}\big(\mathbb{R}; B^{\rho_1}_{r_1, 2}
\big)} &\leq C \|f\|_{H^\mu};\\\label{str2}\big\| K\ast f
\big\|_{L^{q_1}\big(I; B^{\rho_1}_{r_1, 2} \big)} &\leq C\big\|  f
\big\|_{L^{q_2'}\big(I; B^{-\rho_2}_{r'_2, 2} \big)};\\\label{str3}
\big\| K_{R}\ast f \big\|_{L^{q_1}\big(I; B^{\rho_1}_{r_1,  2}
\big)} &\leq C\big\|  f \big\|_{L^{q_2'}\big(I; B^{-\rho_2}_{r'_2,
2} \big)}.
\end{align}
where the subscript $R$ stands for retarded, and
\begin{align*}
K*f&=\int_{\mathbb{R}}K(t-s)f(u(s))ds,\\
K_R\ast f&=\int_{0}^tK(t-s)f(u(s))ds.
\end{align*}

\end{lemma}

\begin{remark}\label{Str-Rem}
 One can check that $(\ref{str2})$ and $(\ref{str3})$ hold  for any
$(q_i,r_i,\theta_i,\rho_i), i=1,2$ satisfying the condition
(\ref{rl}), thus the choice
 of exponents ( especially of $\theta$) is very flexible, which is significant
for the estimate of the nonlinearity. In fact, for any
$(q_1,r_1,\theta_1,\rho_1),(q,r,\theta,\rho)$ satisfying the first
two conditions of (\ref{rl}), we let
$$B:=L^{q_1}\big(I; B^{\rho_1}_{r_1,  2}\big)\cap L^{q}\big(I; B^{\rho}_{r,  2}\big),$$
then  the dual space of $B$ on Hilbert space $H^\mu$ is
 $$B^*:=L^{q'_1}\big(I; B^{2\mu-\rho_1}_{r'_1,  2}\big)\oplus L^{q'}\big(I; B^{2\mu-\rho}_{r',  2}\big).$$
It follows
from $(\ref{str1})$ and the abstract $TT^*$ method that
\begin{equation*}\big\|U\ast f \big\|_{L^{q_1}\big(I; B^{\rho_1}_{r_1,  2}
\big)}\leq\big\| U\ast f \big\|_{B} \leq C\big\| f \big\|_{B^*}\leq
C\big\| f \big\|_{L^{q'}\big(I; B^{2\mu-\rho}_{r', 2} \big)}.
\end{equation*}
where $U\ast f=\int_{\mathbb{R}}U(t-s)f(u(s))ds.$
\end{remark}

\begin{remark}
 According to the above remark, we know that if
$$f(u)=(V*|u|^2)u=\big((V_1 +V_2)*|u|^2\big)u=:f_1(u)+f_2(u),$$
where $V_1\in L^{p_1}, V_2\in L^{p_2}$, then we have
\begin{equation*}
\big\|U*f\big\|_{B} \lesssim \big\| f \big\|_{B^*} \lesssim
\big\|f_1\big\|_{L^{q'_1}\big(I; B^{2\mu-\rho_1}_{r'_1,  2}\big)}+
\big\|f_2\big\|_{L^{q'}\big(I; B^{2\mu-\rho}_{r',  2}\big)}.
\end{equation*}
Hence, without loss of generality, it suffices to consider the
special case  $ V\in L^p$ throughout the paper.
\end{remark}

In addition to the $ST$-norm defined in \eqref{ST}, we also need the
following norm
$$ST^*(I)=[W]^*(I)+[K]^*(I),$$ where
$$[W]^*(I)=L_t^{\frac{2(d+1)}{d+3}}\big(I;B^{\frac12}_{\frac{2(d+1)}{d+3},2}(\mathbb{R}^d)\big),
~[K]^*(I)=L^{\frac{2(d+2)}{d+4}}\big(I;B^{\frac12}_{\frac{2(d+2)}{d+4},2}(\mathbb{R}^d)\big).$$

Now we give a nonlinear estimate which will be applied  to show the
small data scattering which is the first step to obtain the global
time-space estimate that
 lead to the scattering.
\begin{lemma}\label{full}
Let $2<\gamma<\min\{d,4\}$, then we have
\begin{align}\nonumber
&\Big\|\big(|x|^{-\gamma}*|u|^2\big)v\Big\|_{ST^*(I)}+\Big\|\big(|x|^{-\gamma}*(uv)\big)u\Big\|_{ST^*(I)}\\\label{fullst2}
\leq C&\big\| v
\big\|_{[K](I)}\|u\|_{[K](I)}^{\frac4d}\big\|u\big\|_{L_t^\infty(I;
L_x^2)}^{\frac{2(d-2)}{d}}+C\big\| u
\big\|_{[K](I)}^{1+\frac2d}\big\|u\big\|_{L_t^\infty(I;
L_x^2)}^{\frac{d-2}{d}}\big\|v\big\|_{[K](I)}^{\frac2d}\big\|v\big\|_{L_t^\infty(I;
L_x^2)}^{\frac{d-2}{d}}\\\nonumber &+C\big\| v
\big\|_{[W](I)}\|u\|_{L_t^\infty(I;
\dot{H}^{1}_x)}^{\frac{2(d-3)}{d-1}}\|u\|_{[W](I)}
^{\frac{4}{d-1}}+C\big\| u
\big\|_{[W](I)}^{1+\frac{2}{d-1}}\|u\|_{L_t^\infty(I;
\dot{H}_x^1)}^{\frac{d-3}{d-1}}\|v\|_{L_t^\infty(I;
\dot{H}^{1}_x)}^{\frac{d-3}{d-1}}\|v\|_{[W](I)} ^{\frac{2}{d-1}}.
\end{align}
In particular,
\begin{equation}\label{nonlin}
\|(|x|^{-\gamma}*|u|^2)u\|_{ST^*(I)}\leq
C\|u\|_{[k](I)}^{1+\frac4d}\|u\|_{L^\infty(I;H^1)}^{\frac{2(d-2)}{d}}+C\|u\|_{[W](I)}^{1+\frac{4}{d-1}}\|u\|_{L^\infty(I;H^1)}^{\frac{2(d-3)}{d-1}}.
\end{equation}
\end{lemma}

\begin{proof}
We only need to prove the estimate
$\|(|x|^{-\gamma}*|u|^2)v\|_{ST^*(I)}$, since the estimate
$\|(|x|^{-\gamma}*(uv))u\|_{ST^*(I)}$ is similar. From
 the fractional Leibnitz rule and the H\"{o}lder and the
Young inequalities, we have
\begin{align}\nonumber
&\big\|(V*|u|^2)v\big\|_{L^{ q'}\big(I; B^{1/2}_{ {r}', 2}
\big)}\\\label{Bony} \lesssim& \big\| V \big\|_{p} \big\| v
\big\|_{L^{q}\big(I; B^{1/2}_{r, 2} \big)} \big\|u
\big\|^2_{L^k\big(I; L^s \big)}+ \big\| V \big\|_{p} \big\| u
\big\|_{L^{q}\big(I; B^{1/2}_{r, 2} \big)} \big\|u
\big\|_{L^k\big(I; L^s \big)}\big\|v \big\|_{L^k\big(I; L^s \big)},
\end{align}
where we have assumed for simplicity that $V \in L^p$, and the
exponents satisfy
\begin{equation}\label{er}
\left\{ \aligned \frac{d}{p} & = 2 \delta(r) + 2\delta(s), \\
\frac{2}{q}& + \frac{2}{k} =1.
\endaligned \right.
\end{equation}
If we set $V_1=|x|^{-\gamma}\chi_{|x|\leq 1},
V_2=|x|^{-\gamma}\chi_{|x|\geq 1},$ then
\begin{align}\nonumber
&\|(|x|^{-\gamma}*|u|^2)v\|_{ST^*(I)}\\\nonumber
\leq&\|(V_1*|u|^2)v\|_{[W]^*(I)}+\|(V_2*|u|^2)v\|_{[K]^*(I)}\\\label{eqq1.1}
\triangleq&I_1+I_2.
\end{align}
For $I_2$, since $V_2=|x|^{-\gamma}\chi_{|x|\geq 1}\in
L^{\frac{d}{2}}$, if we take admissible pair $q=r=\frac{2(d+2)}{d}$
and $\delta(s_2)=\frac2{k_2}$ (then $\delta(r)=\frac{d}{d+2},
k_2=d+2$) in \eqref{er}, then
\begin{align}\nonumber
I_2&=\|(V_2*|u|^2)v\|_{L^{q^\prime}\big(I;B^{\frac12}_{q^\prime,2}\big)}\\\label{eq1.2}
&\lesssim\big\| v \big\|_{[K](I)} \big\|u \big\|^2_{L^{k_2}(I;
L^{s_2} )}+ \big\| u \big\|_{[K](I)} \big\|u \big\|_{L^{k_2}(I;
L^{s_2} )}\big\|v \big\|_{L^{k_2}(I; L^{s_2} )}.
\end{align}
Using H\"older's inequality and the Sobolev embedding theorem, we
get
\begin{equation}\label{eq1.1}\|v \big\|_{L^{k_2}(I; L^{s_2}
)}\leq\|v\|_{[K](I)}^{\frac2d}\|v\|_{L^\infty(I,
L^2)}^{\frac{d-2}{d}}.
\end{equation}
Plugging \eqref{eq1.1} into \eqref{eq1.2}, we obtain
\begin{equation}\label{eq1.3}I_2\lesssim\big\| v \big\|_{[K](I)}\|u\|_{[K](I)}^{\frac4d}\|u\|_{L_t^\infty
L_x^2}^{\frac{2(d-2)}{d}}+\big\| u
\big\|_{[K](I)}^{1+\frac2d}\|u\|_{L_t^\infty
L_x^2}^{\frac{d-2}{d}}\|v\|_{[K](I)}^{\frac2d}\|v\|_{L_t^\infty
L_x^2}^{\frac{d-2}{d}}.
\end{equation}

On the other hand, for $d\geq 4$, since
$V_1=|x|^{-\gamma}\chi_{|x|\leq 1}\in L^{\frac{d}{4}}$, if we take
admissible pair $q=r=\frac{2(d+1)}{d-1}$ and
$\delta(s_1)=1+\frac1{k_1}$ (then
$\delta(r)=\frac{d}{d+1},k_1=d+1$), then
\begin{align}\nonumber
I_1&=\|(V_1*|u|^2)v\|_{L^{q^\prime}(I;B^{\frac12}_{q^\prime,2})}\\\label{eq1.4}
&\lesssim\big\| v \big\|_{[W](I)} \big\|u \big\|^2_{L^{k_1}(I;
L^{s_1} )}+ \big\| u \big\|_{[W](I)} \big\|u \big\|_{L^{k_1}(I;
L^{s_1} )}\big\|v \big\|_{L^{k_1}(I; L^{s_1} )}.
\end{align}
The H\"older inequality and the Sobolev embedding theorem yield that
\begin{equation}\label{eq1.5}\|v \big\|_{L^{k_1}(I; L^{s_1}
)}\leq\|v\|_{L_t^\infty
L^{2^*}_x}^{\frac{d-3}{d-1}}\|v\|_{L_t^{\frac{2(d+1)}{d-1}}L_x^{\frac{2d(d+1)}{d^2-2d-1}}}
^{\frac{2}{d-1}}\lesssim\|v\|_{L_t^\infty
\dot{H}^{1}_x}^{\frac{d-3}{d-1}}\|v\|_{[W](I)} ^{\frac{2}{d-1}}.
\end{equation}
Plugging \eqref{eq1.5} into \eqref{eq1.4}, we get
\begin{equation}\label{eq1.6}I_1\lesssim\big\| v \big\|_{[W](I)}\|u\|_{L_t^\infty
\dot{H}^{1}_x}^{\frac{2(d-3)}{d-1}}\|u\|_{[W](I)}
^{\frac{4}{d-1}}+\big\| u
\big\|_{[W](I)}^{1+\frac{2}{d-1}}\|u\|_{L_t^\infty
\dot{H}_x^1}^{\frac{d-3}{d-1}}\|v\|_{L_t^\infty
\dot{H}^{1}_x}^{\frac{d-3}{d-1}}\|v\|_{[W](I)} ^{\frac{2}{d-1}}.
\end{equation}
For $d=3$, since $V_1=|x|^{-\gamma}\chi_{|x|\leq 1}\in L^{1}$, if we
take admissible pair $q=r=k_1=4$ and
$\delta(r)=\delta(s_1)=\frac34$, then
\begin{align}\nonumber
I_1&=\|(V_1*|u|^2)v\|_{L^{\frac43}\big(I;B^{\frac12}_{\frac43,2}\big)}\\\nonumber
&\lesssim\big\| v \big\|_{[W](I)} \big\|u \big\|^2_{L^4(I; L^4 )}+
\big\| u \big\|_{[W](I)} \big\|u \big\|_{L^4(I; L^4 )}\big\|v
\big\|_{L^4(I; L^4 )}\\\label{eq1.7} &\lesssim\big\| v
\big\|_{[W](I)}\big\| u \big\|_{[W](I)}^2.
\end{align}
Combining \eqref{eqq1.1}, \eqref{eq1.3}, \eqref{eq1.6} and
\eqref{eq1.7}, we obtain
\begin{align}\nonumber
&\big\|(|x|^{-\gamma}*|u|^2)v\big\|_{ST^*(I)}\\\nonumber
\lesssim&\big\| v
\big\|_{[K](I)}\|u\|_{[K](I)}^{\frac4d}\|u\|_{L_t^\infty
L_x^2}^{\frac{2(d-2)}{d}}+\big\| u
\big\|_{[K](I)}^{1+\frac2d}\|u\|_{L_t^\infty
L_x^2}^{\frac{d-2}{d}}\|v\|_{[K](I)}^{\frac2d}\|v\|_{L_t^\infty
L_x^2}^{\frac{d-2}{d}}\\\label{eq1.8} &+\big\| v
\big\|_{[W](I)}\|u\|_{L_t^\infty
\dot{H}^{1}_x}^{\frac{2(d-3)}{d-1}}\|u\|_{[W](I)}
^{\frac{4}{d-1}}+\big\| u
\big\|_{[W](I)}^{1+\frac{2}{d-1}}\|u\|_{L_t^\infty
\dot{H}_x^1}^{\frac{d-3}{d-1}}\|v\|_{L_t^\infty
\dot{H}^{1}_x}^{\frac{d-3}{d-1}}\|v\|_{[W](I)} ^{\frac{2}{d-1}},
\end{align}
which completes the proof of Lemma \ref{full}.
\end{proof}

We can now state the local well-posedness for $(\ref{equ1})$ with
large initial data and small data scattering in the energy space.

\begin{theorem}\label{small}
Assume $d\geq 3,~2<\gamma<\min\{d,4\}$  and $(u_0,u_1)\in
H^1(\mathbb{R}^d)\times L^2(\mathbb{R}^d)$. There exists a small
constant $\delta=\delta(E)$ such that if $\|(u_0,u_1)\|_{
H^1\times L^2}\leq E$ and $I$ is an interval such that
$$\|\dot{K}(t)u_0 + K(t)u_1\|_{ST(I)}\leq
\delta,$$ then there exists a unique strong solution $u$ to
\eqref{equ1} in $I\times \mathbb{R}^d$, with $u\in C(I;H^1)\cap
C^1(I;L^2)$ and
\begin{equation}\label{smalll}
\|u\|_{ST(I)}\leq 2C\delta. \end{equation} Let $(-T_{min},T_{max})$
be the maximal time interval on which $u$ is well-defined. Then, if
$T_{max}<+\infty,$ then $\lim\limits_{t\rightarrow
T_{max}}\|u(t)\|_{H^1}=+\infty.$ Similarly, if\ $T_{min}<+\infty,$
then $\lim\limits_{t\rightarrow -T_{min}}\|u(t)\|_{H^1}=+\infty.$
\end{theorem}
\begin{proof}
We apply the Banach fixed point argument to prove this lemma. First
we define the map
\begin{equation}\label{inte3}
\Phi(u(t))=\dot{K}(t)u_0 + K(t)u_1-\int^{t}_{0}K(t-s)f(u(s))ds
\end{equation}
on the complete metric space $B$
\begin{align*}
B=\{u\in C(I;H^1):\|u\|_{L^\infty(I;H^1)}\leq
2C\|(u_0,u_1)\|_{H^1\times L^2},\|u\|_{ST(I)}\leq2C\delta\}
\end{align*}
with the metric $d(u,v)=\big\|u-v\big\|_{ST(I)\cap
L^\infty(I;H^1)}$.

It suffices to prove that the operator defined by the RHS of
$(\ref{inte3})$ is a contraction map on $B$ for $I$. In fact, if
$u\in B$, then  by Lemma $\ref{lem22}$ and \eqref{nonlin}, we have
\begin{align*}\|\Phi(u)\|_{ST(I)}\leq& C\|\dot{K}(t)u_0 + K(t)u_1\|_{ST(I)}+C\|f(u)\|_{ST^*(I)}\\
\leq&
C\delta+C\|u\|_{[k](I)}^{1+\frac4d}\|u\|_{L^\infty(I;H^1)}^{\frac{2(d-2)}{d}}+C\|u\|_{[W](I)}
^{1+\frac{4}{d-1}}\|u\|_{L^\infty(I;H^1)}^{\frac{2(d-3)}{d-1}}\\
\leq &C\delta+C(2\delta)^{\frac{d+4}{d}}(2C\|(u_0,u_1)\|_{H^1\times
L^2})^{\frac{2(d-2)}{d}}\\
&+C(2\delta)^{\frac{d+3}{d-1}}(2C\|(u_0,u_1)\|_{H^1\times
L^2})^{\frac{2(d-3)}{d-1}}.
\end{align*}
Hence, if we take $\delta$ sufficiently small such that
$C2^{\frac{d+4}{d}}\delta^\frac{4}{d}(2C\|(u_0,u_1)\|_{H^1\times
L^2})^{\frac{2(d-2)}{d}}\leq\frac12$ and
$C(2\delta)^{\frac{4}{d-1}}(2C\|(u_0,u_1)\|_{H^1\times
L^2})^{\frac{2(d-3)}{d-1}}\leq\frac12$, then
$\|\Phi(u)\|_{ST(I)}\leq 2C\delta.$

Similarly, we get $\|\Phi(u)\|_{L^\infty(I;H^1)}\leq
2C\|(u_0,u_1)\|_{H^1\times L^2}$, and so $\Phi(u)\in B.$

On the other hand,  for $\omega_1, \omega_2\in B$, by Strichartz
estimate, we obtain
\begin{equation*}\begin{split}
d(\Phi(\omega_{1}),\Phi(\omega_2))\leq&
C\big\|(|x|^{-\gamma}*|\omega_1|^2)\omega_1-(|x|^{-\gamma}*|\omega_2|^2)\omega_2\big\|_{ST^*(I)}\\
\leq&
C\big\|(|x|^{-\gamma}*|\omega_1|^2)(\omega_1-\omega_2)\big\|_{ST^*(I)}+C\big\|\big[|x|^{-\gamma}*(\omega_1-\omega_2)^2\big]
\omega_2\big\|_{ST^*(I)}\\
&+2C\big\|\big[|x|^{-\gamma}*(\omega_2(\omega_1-\omega_2))\big]
\omega_2\big\|_{ST^*(I)}.
\end{split}
\end{equation*}
By Lemma \ref{full}, the above quantity can be controlled by
\begin{align*}
Cd(\omega_1,\omega_2)\sum\limits_{i,j=1}^2\Big(&\|\omega_i\|_{ST(I)}^\frac4d\|\omega_j\|_{L^\infty_t
L_x^2}^\frac{2(d-2)}{d}+\|\omega_i\|_{ST(I)}^\frac{d+2}{d}\|\omega_j\|_{L^\infty_t
L_x^2}^\frac{(d-2)}{d}\\
&+\|\omega_i\|_{ST(I)}^\frac4{d-1}\|\omega_j\|_{L^\infty_t
H_x^1}^\frac{2(d-3)}{d-1}+\|\omega_i\|_{ST(I)}^\frac{d+1}{d-1}\|\omega_j\|_{L^\infty_t
H_x^1}^\frac{(d-3)}{d-1}\Big),
\end{align*}
which allows us to derive
\begin{equation*}
d(\Phi(\omega_{1}),\Phi(\omega_2))\leq\frac{1}{2}d(\omega_1,\omega_2)
\end{equation*}
by taking $\delta$ small. This completes the proof.
\end{proof}
In the defocusing case ($\mu=1$), using the conservation of energy,
the above solutions can be extended globally. This gives the
following proposition which is the starting point of our
investigation. As we will see in Corollary \ref{global}, the above
solutions can also be extended globally for the focusing case
($\mu=-1$) under the restriction $E(u_0,u_1)<E(\mathcal{W},0)$ and $\|\nabla
u_0\|_2^2+\|u_0\|_2^2<\|\nabla \mathcal{W}\|_2^2+\|\mathcal{W}\|_2^2$.

\begin{proposition}\label{pro23}
Assume $d\geq 3$, $\mu=1$ and $2<\gamma<\min\{4,d\}$. For all
initial data $(u_0, u_1)\in H^1 \times L^2$, there exists a unique
globally defined nonlinear solution $u\in C(\mathbb{R}, H^1)\cap
C^1(\mathbb{R}, L^2)\cap ST_{loc}[\mathbb{R}]$. Besides, the
evolution flow $(u(0),u_t(0))\in H^1 \times L^2\mapsto u$ is
continuous for all compact time interval $I$.
\end{proposition}
The above Proposition only gives the local in time bounds, while the
key point to understand better the behaviour of these solutions is
to gain access to global in time bounds.  Actually, as mentioned in
the introduction, a global in time bound of the $ST$- norm of the
solution is sufficient for scattering, this is the object of the
following proposition.

\begin{proposition}\label{prop2.3}
Let $d\geq 3$, $\mu=1$, $2<\gamma<\min\{4,d\}$ and $u\in
C(\mathbb{R},H^1(\mathbb{R}^d))\cap
C^1(\mathbb{R},L^2(\mathbb{R}^d))$ be a global strong solution of
\eqref{equ1} with initial data $(u_0,u_1)\in H^1\times L^2$.  Then
we have if
$$\|u\|_{ST(\mathbb{R})}<\infty,$$ then $u$ scatters.

Moreover, in the focusing case: $\mu=-1$, the above result also
holds under the restriction $E(u_0,u_1)<E(\mathcal{W},0),$ and $\|\nabla
u_0\|_2^2+\|u_0\|_2^2<\|\nabla \mathcal{W}\|_2^2+\|\mathcal{W}\|_2^2$.
\end{proposition}

\begin{proof}
We just prove that $u$ scatters at $+\infty$, the proof for the
scattering at $-\infty$ is similar.  Using Duhaml's formula, the
solution with initial data $(u(0),\dot{u}(0))=(u_0,u_1)\in H^1\times
L^2$ of \eqref{equ1} can be written as
\begin{equation}
{u(t)\choose\dot{u}(t)}=V_0(t){u_0 \choose
u_1}-\int^{t}_{0}V_0(t-s){0\choose f(u(s))}ds,
\end{equation}
where $V_0$ is defined by \eqref{inte2}. Denote the scattering data
$(u_0^+,u_1^+)$ by
$${u_0^+\choose u_1^+}={u_0\choose u_1}-\int_0^\infty V_0(-s){0\choose\mu(|x|^{-\gamma}*|u|^2)u(s)}ds.$$
Then, by Strichartz estimate \eqref{str3} and \eqref{eq1.8},  we can
obtain
\begin{align*}
\Big\|{u\choose \dot{u}}-V_0(t){u_0^+\choose u_1^+}\Big\|_{H^1\times L^2}&=\Big\|\int_t^\infty V_0(t-s){0\choose f(u(s))}ds\Big\|_{H^1\times L^2}\\
&\lesssim\|(|x|^{-\gamma}*|u|^2)u\|_{ST^*(t,\infty)}\\
&\lesssim\|u\|_{[K](t,\infty)}^{1+\frac{4}{d}}\|u\|_{L_t^\infty
L_x^2}^{\frac{2(d-2)}{d}}+\|u\|_{[W](t,\infty)}^{1+\frac{4}{d-1}}\|u\|_{L_t^\infty
H_x^1}^{\frac{2(d-3)}{d-1}}\\
&\rightarrow 0,\quad \text{as}\quad t\rightarrow\infty.
\end{align*}
Thus $u$ scatters. Here we used the fact that: $\|u\|_{L_t^\infty
(\mathbb{R}; H^1)}$ can be controlled by $E(u_0,u_1)$ from the
conservation of the energy in the defocusing case, and
$\|u\|_{L_t^\infty (\mathbb{R}; H^1)}$ can also be controlled by
$E(u_0,u_1)$ for the focusing case by Lemma \ref{freeequi} below .
\end{proof}

\subsection{Perturbation lemma}
 We here record the short and long time perturbations as in
\cite{CKSTT07}. Roughly speaking, the stability proposition says
that, if initial data are close enough and the perturbation term is
small in some sense, then the solutions will be close.

\begin{lemma}[Short-time perturbations]\label{short}  Let $I$ be a
time interval, and let $\tilde{u}$ be a function on $I\times \R^d$
which is  a near solution to (\ref{equ1}) in the sense that
\begin{equation}\label{near}
(\square+1)\tilde{u}=-f(\tilde{u})+e
\end{equation}
for some function $e$. Assume that
\begin{eqnarray*}
 \|\tilde{u}\|_{L_t^\infty
H^1_x(I\times \R^d)}+\|\partial_{t}\tilde{u}\|_{L_t^\infty
L^2_x(I\times \R^d)} \leq E
\end{eqnarray*}
for some constant $E>0$. Let $t_0\in I$, and let $(u(t_0),u_t(t_0))$
be close to $(\tilde{u}(t_0),\tilde{u}_t(t_0))$ in the sense that
\begin{equation}\label{eq2.14}
\|(u(t_0)-\tilde{u}(t_0),u_t(t_0)-\tilde{u}_t(t_0))\|_{H^1\times
L^2}\leq E'
\end{equation}
and assume also that  we have smallness conditions
\begin{align}\label{eq2.15}
\|\dot{K}(t-t_0)(u(t_0)-\tilde{u}(t_0))+K(t-t_0)(u_t(t_0)-\tilde{u}_t(t_0))\|_{ST(I)}\leq\epsilon,\\\label{equ2.16}
\|\tilde{u}\|_{ST(I)}\leq \epsilon, \quad \|e\|_{ST^*(I)}\leq
\epsilon
\end{align}
for some $0<\epsilon<\epsilon_0$, where $\epsilon_0=\epsilon_0(E)>0$
is a small enough constant.

We conclude that there exists a solution $u$ to (\ref{equ1}) on
$I\times\R^d$ with the specified initial data $(u(t_0),u_t(t_0))$ at
$t_0$, and furthermore
\begin{equation} \label{eq2.16}
 \begin{aligned}
\|u-\tilde{u}\|_{ST(I)}\lesssim \epsilon,\\
\|(\square+1)(u-\tilde{u})\|_{ST^*(I)}\lesssim \epsilon,
\end{aligned}
\end{equation} where $\square:=\partial_{tt}-\Delta$.
\end{lemma}

\begin{proof}First we claim that
$$\|\tilde{u}\|_{L^{k_i}(I;L^{s_i})}\leq \epsilon_i,~(i=1,2),$$
where $\epsilon_i$ is a small constant depending on $\epsilon_0$ and
$E$, $k_i,s_i$ defined in \eqref{eq1.2}, \eqref{eq1.4} and
\eqref{eq1.7}. In fact, by \eqref{eq1.1} and \eqref{eq1.5}, we
obtain
\begin{align*}
\|\tilde{u}\|_{L^{k_i}(I;L^{s_i})}&\leq
\|\tilde{u}\|_{ST(I)}^{\theta_i}\|\tilde{u}\|_{L^\infty(I;H^1)}^{1-\theta_i}\\
&\lesssim \epsilon_0^{\theta_i} E^{1-\theta_i}\triangleq\epsilon_i.
\end{align*}
This concludes the claim. Similarly,
$$\|\dot{K}(t-t_0)(u(t_0)-\tilde{u}(t_0))+K(t-t_0)(u_t(t_0)-\tilde{u}_t(t_0))\|_{L^{k_i}(I;L^{s_i})}\leq \epsilon_i,~(i=1,2).$$

 Let
$w:=u-\tilde {u}$ and $Y(I)=ST(I)\cap L^{k_1}(I;L^{s_1})\cap
L^{k_2}(I;L^{s_2})$, then by Lemma \ref{lem22}, we have
\begin{equation}\label{eq2.17}
 \|w\|_{Y(I)}\lesssim \epsilon+\|(\square+1)w\|_{ST^*(I)}=:\epsilon+S(I),\end{equation}
 where
\begin{equation}\label{eq2.18}
(\square+1)w=-f(u)+f(\tilde {u})-e.
\end{equation}
Hence, by Lemma \ref{full} and \eqref{Bony}, we have
\begin{equation*}\begin{split}
S(I)\leq&\big\|f(u)-f(\tilde {u})\big\|_{ST^*(I)}+\big\|e\big\|_{ST^*(I)}\\
\lesssim&
\big\|(|x|^{-\gamma}*|\tilde{u}|^2)w\big\|_{ST^*(I)}+\big\|(|x|^{-\gamma}*|w|^2)\tilde{u}\big\|_{ST^*(I)}+\big\|(|x|^{-\gamma}*|w|^2)w\big\|_{ST^*(I)}\\
&+\big\|(|x|^{-\gamma}*(w\tilde{u}))w\big\|_{ST^*(I)}+\big\|(|x|^{-\gamma}*(w\tilde{u}))\tilde{u}\big\|_{ST^*(I)}+\epsilon\\
\lesssim&\sum_{l=1}^3\|w\|_{Y(I)}^l\|\tilde{u}\|_{Y(I)}^{3-l}+\epsilon\\
\lesssim& \sum_{l=1}^3 (S(I)+\epsilon)^l
(\epsilon+\epsilon_1+\epsilon_2)^{3-l}+\epsilon.\end{split}
\end{equation*}

Making use of a standard continuity argument, it follows that
$S(I)\lesssim\epsilon$ and then $\|w\|_{Y(I)}\lesssim \epsilon$.
This, together with \eqref{equ2.16},
 implies the first estimate in  \eqref{eq2.16}. This completes the proof.
\end{proof}

\begin{lemma}[Long-time perturbations]\label{long}  Let $I$ be a
time interval, and let $\tilde{u}$ be a function on $I\times \R^n$
which is  a solution to \eqref{near}  such that
\begin{equation}\label{eq2.20}\begin{aligned}
\|\tilde{u}\|_{ST(I)}\leq M,\\ \|\tilde{u}\|_{L_t^\infty
H^1_x(I\times \R^d)}+\|\partial_{t}\tilde{u}\|_{L_t^\infty
L^2_x(I\times \R^d)}\leq E,
\end{aligned}
\end{equation}
for some constant $M,E>0$. Let $t_0\in I$, and let
$(u(t_0),u_t(t_0))$ be close to $(\tilde{u}(t_0),\tilde{u_t}(t_0))$
in the sense that
\begin{equation}\label{eq2.21}
\big\|\big(u(t_0)-\tilde{u}(t_0),u_t(t_0)-\tilde{u}_t(t_0)\big)\big\|_{H^1\times
L^2}\leq E',
\end{equation}and assume also that  we have smallness conditions
\begin{equation} \label{eq2.22}
\big\|\dot{K}(t-t_0)(u(t_0)-\tilde{u}(t_0))+K(t-t_0)(u_t(t_0)-\tilde{u}_t(t_0))\big\|_{ST(I)}+\|e\|_{ST^*(I)}\leq\epsilon
\end{equation}
for some small $0<\epsilon<\epsilon_1$, where
$\epsilon_1=\epsilon_1( M, E)>0$.  We conclude that there exists a
solution  $u(t)$ to \eqref{equ1} on $I\times\R^d$ with specific
initial data $(u(t_0),u_t(t_0))$ at $t_0$, and furthermore
\begin{equation}\label{eq2.23}
\begin{aligned}
\|u-\tilde{u}\|_{ST(I)}\leq & C(M,E)\epsilon,\\
\|u\|_{ST(I)}\leq & C(M,E,E').
\end{aligned}
\end{equation}
\end{lemma}

\begin{proof} Since $\|\tilde{u}\|_{ST(I)}\leq M$, we may subdivide $I$
into $N=C(M, \varepsilon_0)$  intervals $I_j=[t_j,t_{j+1}]$ such
that
$$\|\tilde{u}\|_{ST(I_j)}\leq \epsilon_0, \quad \quad 1\le j\le C(M,\epsilon_0), $$
where $\epsilon_0=\epsilon_0(E)>0$ is the same as Lemma \ref{short}.

Next we can use inductively the short-time perturbations lemma for
$j=0,1,\cdots,N$ to get
\begin{align*}
\|u-\tilde{u}\|_{ST(I_j)}&\leq C(j) \epsilon,\\
\|(\square+1)(u-\tilde{u})\|_{ST^*(I_j)}&\leq C(j) \epsilon
\end{align*}
and then we obtain
\begin{equation} \label{eq2.27}
\begin{split}
&\big\|(u(t_{j+1})-\tilde{u}(t_{j+1}),\dot{u}(t_{j+1})-\dot{\tilde{u}}(t_{j+1}))\big\|_{H^1\times L^2}\\
\leq
&\big\|(u(t_{j})-\tilde{u}(t_{j}),\dot{u}(t_{j})-\dot{\tilde{u}}(t_{j}))\big\|_{H^1\times
L^2}+C(j) \epsilon,
\end{split}
\end{equation}
thus the claim follows by the standard argument.
\end{proof}

\section{Variational characterizations}
\setcounter{section}{3}\setcounter{equation}{0} In this section, we
discuss some properties of ground states,
and some preliminary lemmas for the study in the focusing case. The
idea is similar to S. Ibrahim, N. Masmoudi, K. Nakanishi \cite{IMN} and
C. Miao, Y. Wu \cite{MW}.

First, by  the symmetry
$$\iint_{\mathbb{R}^d\times\mathbb{R}^d}\frac{x\cdot(x-y)}{|x-y|^{\gamma+2}}\phi(x)^2\phi(y)^2dxdy
=-\iint_{\mathbb{R}^d\times\mathbb{R}^d}\frac{y\cdot(x-y)}{|x-y|^{\gamma+2}}\phi(x)^2\phi(y)^2dxdy,$$
and a direct computation we have the following identities:
\begin{lemma}\label{Pohozaev}
Assume $\phi\in\mathcal{S}(\mathbb{R}^d)$, then
\begin{align*}
-\int_{\mathbb{R}^d} \Delta\phi
x\cdot\nabla\phi dx&=-\frac{d-2}{2}\int_{\mathbb{R}^d}|\nabla\phi|^2dx,\\
\int_{\mathbb{R}^d}\phi x\cdot\nabla\phi dx&=-\frac{d}{2}\int_{\mathbb{R}^d}|\phi|^2dx,\\
\int_{\mathbb{R}^d}(|x|^{-\gamma}*|\phi|^2)\phi x\cdot\nabla\phi
dx&=\big(-\frac{d}{2}+\frac{\gamma}{4}\big)\iint_{\mathbb{R}^d\times
\mathbb{R}^d}\frac{\phi(x)^2\phi(y)^2}{|x-y|^{\gamma}}dxdy.
\end{align*}
\end{lemma}

\begin{lemma}\label{lemma3.1}
Assume that $2<\gamma<\min\{4,d\}$. Let $\phi$ be the
$H^1(\mathbb{R}^d)$ solution of the following equation
\begin{equation}\label{ground}
-\Delta \phi+\phi=(|x|^{-\gamma}*|\phi|^2)\phi,
\end{equation}
then the following identity holds:
$$K_1(\phi)\triangleq\int_{\mathbb{R}^d}\big(|\nabla\phi|^2+|\phi|^2\big)dx-\iint_{\mathbb{R}^d\times
\mathbb{R}^d}\frac{\phi(x)^2\phi(y)^2}{|x-y|^\gamma}dxdy=0.$$
$$K_2(\phi)\triangleq\int_{\mathbb{R}^d}|\nabla\phi|^2dx-\frac{\gamma}{4}\iint_{\mathbb{R}^d\times
\mathbb{R}^d}\frac{\phi(x)^2\phi(y)^2}{|x-y|^{\gamma}}dxdy=0.$$
\end{lemma}
\begin{proof}
$K_1(\phi)=0$ can be obtained by multiplying \eqref{ground} both
sides by $\phi$ and integrating.  By Lemma \ref{Pohozaev},
$K_2(\phi)=0$ is obtained by multiplying \eqref{ground} both sides
by $x\cdot\nabla\phi+\frac{d}{2}\phi$, and integrating.

\end{proof}
Let the static energy $J$ be defined by
\begin{equation}\label{static}
J(\phi)=\frac12\int_{\mathbb{R}^d}\big[|\nabla\phi|^2+|\phi|^2\big]dx-\frac14\iint_{\mathbb{R}^d\times
\mathbb{R}^d}\frac{\phi(x)^2\phi(y)^2}{|x-y|^\gamma}dxdy.
\end{equation}
Let $\Omega\triangleq\{\varphi\in H^1(\mathbb{R}^d)\backslash\{0\}:
\varphi\ \text{solves}\ \eqref{ground}\}$, then the ground state set
of the elliptic equation \eqref{ground} is defined as
$$\Lambda\triangleq\{\varphi\in\Omega:\ J(\varphi)\leq J(\phi),\
\forall \phi\in\Omega\}.$$ The existence of the ground state was
shown in \cite{Liu} where it had been shown that the ground state a
radial, rapidly decaying function. And let $\mathcal{W}$ be an element of the
ground state set. For our purpose, we will give two characterizations
of the ground state based on the functional $K_1, K_2$, which will
be important to describe the structures of the dichotomy of blow
up and scattering associated to the nonlinear Klein-Gordon-Hartree
equation.

\begin{lemma}\label{noempty}
Let $\Omega_1=\{\phi\in H^1(\mathbb{R}^d):\
\phi\not=0,K_1(\phi)=0\},$ and
$$\Lambda_1\triangleq\{\varphi\in\Omega_1:\ J(\varphi)\leq J(\phi),\
\forall \phi\in\Omega_1\},\ $$
\begin{equation}\label{mini1}m\triangleq\inf\{J(\phi):\phi\in H^1(\mathbb{R}^d),\
\phi\not=0,\ K_1(\phi)= 0\}\end{equation} then
$\Lambda_1\neq\emptyset,$ and moreover, $m=J(\mathcal{W}),$ i.e. $m$ is
attained by the ground state $\mathcal{W}$.
\end{lemma}
Before proving Lemma \ref{noempty}, we introduce some notations.
First we decompose $K_j(\phi)\ (j=1,2)$ into the quadratic and
nonlinear parts:
\begin{align}
\begin{cases}
K_j(\phi)=K_j^Q(\phi)+K_j^N(\phi),\ j=1,2,\\
K_1^Q(\phi)=\|\phi\|^2_{L^2}+\|\phi\|^2_{\dot{H}^1(\mathbb{R}^d)},\ K_2^Q(\phi)=\|\nabla\phi\|_2^2,\\
K_1^N(\phi)=-\iint_{\mathbb{R}^d\times
\mathbb{R}^d}\frac{\phi(x)^2\phi(y)^2}{|x-y|^\gamma}dxdy,\
K_2^N(\phi)=\frac{\gamma}{4}K_1^N(\phi).
\end{cases}\end{align}
\begin{remark}\label{infty10} If $\phi\in H^1(\mathbb{R}^d),$ then
$K_1^Q(e^\lambda\phi)=e^{2\lambda}\big(\|\phi\|^2_{L^2}+\|\phi\|^2_{\dot{H}^1}\big)\rightarrow
0,$ as $\lambda\rightarrow-\infty.$
\end{remark}
\begin{lemma}\label{positive}
Assume $2<\gamma<\min\{4,d\},$ then for any bounded sequence
$\{\phi_n\}\subset H^1(\mathbb{R}^d)\backslash\{0\}$ such that
$K_1^Q(\phi_n)\rightarrow 0$, we have for large $n$,
$$K_1(\phi_n)>0.$$
\end{lemma}
\begin{proof}
Using the H\"older and generalized Young inequality, we have
\begin{equation}\label{gamma1}
|K_1^N(\phi)|=\big\|(|x|^{-\gamma}*|\phi|^2)\phi^2\big\|_{L^1}
\leq\|\phi\|_{L^{2p}}^4\big\||x|^{-\gamma}\big\|_{L^{\frac{d}{\gamma}}_*},
\end{equation}
here
$p=\frac{2d}{2d-\gamma},1+\frac{1}{p^\prime}=\frac{\gamma}{d}+\frac1p,$
and by the H\"older inequality and Sobolev imbedding theorem, one
has
\begin{equation}\label{gamma}
\|\phi\|_{L^{2p}}\leq\|\phi\|_{L^2}^{1-\frac{\gamma}{4}}\|\phi\|_{\frac{2d}{d-2}}^{\frac{\gamma}{4}}\lesssim\|\phi\|_{L^2}^{\frac{4-\gamma}{4}}
\|\phi\|_{\dot{H}^1}^{\frac{\gamma}{4}}.\end{equation} Plugging
\eqref{gamma} into \eqref{gamma1}, we obtain
\begin{equation}\label{gamma2}
|K_1^N(\phi)|\lesssim\|\phi\|_{L^2}^{4-\gamma}\|\phi\|_{\dot{H}^1}^{\gamma}\lesssim
K_1^Q(\phi)^2.
\end{equation}
Since $\phi_n\in H^1(\mathbb{R}^d)\backslash\{0\}$ satisfies
$K_1^Q(\phi_n)\rightarrow 0$, \eqref{gamma2} and
$K_1(\phi_n)=K_1^Q(\phi_n)+K_1^N(\phi_n)$, we get for large $n$,
$$K_1(\phi_n)>0.$$
\end{proof}

\begin{lemma}\label{positive1}
If we set $H_1(\phi)\triangleq-\frac14
K_1^N(\phi)=\frac14\int\int_{\mathbb{R}^d\times\mathbb{R}^d}\frac{\phi(x)^2\phi(y)^2}{|x-y|^\gamma}dxdy,$
and \begin{equation}\label{mini} m_1=\inf\{H_1(\phi):\ \phi\in
H^1(\mathbb{R}^d),\ \phi\not=0,\ K_1(\phi)\leq 0\},
\end{equation}
then $m_1=m$.
\end{lemma}

\begin{proof}
It is trivial to prove that $m_1\leq m$ because $J(\phi)=H_1(\phi)$
if $K_1(\phi)=0.$ So it suffices to show $m\leq m_1$. Take $\phi\in
H^1$ such that $K_1(\phi)<0$.

It follows from Remark \ref{infty10} and Lemma \ref{positive} that
$$\exists \lambda_1<0,\ s.t.\quad K_1(e^{\lambda_1}\phi)>0.$$
Combining this with $K_1(\phi)<0$, we deduce that
\begin{equation}\label{pos1}\exists \lambda\in(\lambda_1,0),\ s.t.\ K_1(e^{\lambda}\phi)=0,\
J(e^\lambda\phi)=H_1(e^\lambda\phi)=e^{4\lambda}H_1(\phi)\leq
H_1(\phi),
\end{equation}
and so  $m\leq m_1.$ This completes the proof of Lemma
\ref{positive1}.
\end{proof}
\noindent{\bf The proof of Lemma \ref{noempty}:}

 {\bf Step1:
} We claim that $\Lambda_1\neq\emptyset.$ Indeed, let $\phi_n\in
H^1$ be a minimizing sequence for \eqref{mini}, i.e.
$$K_1(\phi_n)\leq 0,\ \phi_n\not=0,\ H_1(\phi_n)\searrow m.$$
Let $\phi_n^*$ be the Schwartz symmetrization of $\phi_n$, i.e. the
radial decreasing rearrangement(see \cite{Lieb}). Since
$\|\nabla\phi_n^*\|_2\leq\|\nabla\phi_n\|_2$,
$\|\phi_n^*\|_2=\|\phi_n\|_2$ and
$$\iint_{\mathbb{R}^d\times \mathbb{R}^d}\frac{\phi_n(x)^2\phi_n(y)^2}{|x-y|^\gamma}dxdy\leq\iint_{\mathbb{R}^d\times \mathbb{R}^d}
\frac{|\phi_n^*(x)|^2|\phi_n^*(y)|^2}{|x-y|^\gamma}dxdy,$$ (by the
General rearrangement inequality, see page 93 Theorem3.8 in
\cite{Lieb}), therefore we have
$$\phi_n^*\not=0,\ K_1(\phi_n^*)\leq K_1(\phi_n)\leq 0,\ H_1(\phi_n^*)\geq H_1(\phi_n).$$
If we choose $\lambda_n$ such that
$e^{4\lambda_n}=\frac{H_1(\phi_n)}{H_1(\phi_n^*)}$, then
$$\phi_n^*\not=0,\ K_1(e^{\lambda_n}\phi_n^*)\leq 0,\
H_1(\phi_n)=H_1(e^{\lambda_n}\phi_n^*)\searrow m.$$Then by
\eqref{pos1}, we may replace it by symmetric $\varphi_n\in H^1$ such
that
\begin{equation}\label{lim}
\varphi_n\not=0,\ K_1(\varphi_n)=0,\
J(\varphi_n)=H_1(\varphi_n)\rightarrow m.
\end{equation}
Notice that
$J(\varphi_n)=\frac14(\|\nabla\varphi_n\|_2^2+\|\varphi_n\|_2^2)\rightarrow
m,$ we know that $\{\varphi_n\}$ is bounded in $H^1$. And so up to
subsequence, it converges to some $\varphi$ weakly in $H^1$. By the
radial symmetry, it also converges strongly in $L^p$ for all
$2<p<2^*=\frac{2d}{d-2}$. By the  H\"older and Young inequality, we
deduce that
$$\big\|(|x|^{-\gamma}*|\varphi_n|^2)\varphi_n^2-(|x|^{-\gamma}*|\varphi|^2)\varphi^2\big\|_{L^1}
\lesssim\big\|\varphi_n-\varphi\big\|_p\big\|\varphi_n+\varphi\big\|_p
\Big(\big\|\varphi_n\big\|_p^2+\big\|\varphi\big\|_p^2\Big),$$ where
$2<p=\frac{4d}{2d-\gamma}<2^*$, and so the nonlinear parts
$K_1^N(\varphi_n)$ converges. And by the Fatou Lemma, we have
$$K_1(\varphi)\leq\varliminf\limits_{n\rightarrow\infty}K_1(\varphi_n)\leq
0,\ \text{and}\
H_1(\varphi)\leq\varliminf\limits_{n\rightarrow\infty}H_1(\varphi_n)\leq
m.$$ If $\varphi=0$, then $K_1(\varphi_n)=0$ implies that
$K_1^Q(\varphi_n)=-K_1^N(\varphi_n)\rightarrow -K_1^N(\varphi)=0$ by
$\varphi=0$, and by Lemma \ref{positive} we have $K_1(\varphi_n)>0$
for large $n$, which contradicts with $K_1(\varphi_n)=0$. Hence
$\varphi\not=0.$

By \eqref{pos1}, we may replace $\varphi$ by $e^\lambda\varphi$, so
that $K_1(\varphi)=0, J(\varphi)=H_1(\varphi)\leq m$ and
$\varphi\not=0$. Then $\varphi$ is a minimizer and
$m=H_1(\varphi)>0.$ Hence $\Lambda_1\not=\emptyset$.

{\bf Step2: $m=J(\mathcal{W})$.} That is, we need to prove $J(\varphi)=J(\mathcal{W})$,
where $\varphi$ is attained in Step 1. Since $\varphi$ is a minimizer
for \eqref{mini1}, there exists a Lagrange multiplier
$\eta\in\mathbb{R}$ such that
\begin{equation}\label{lan1}
J^\prime(\varphi)=\eta K_1^\prime(\varphi).
\end{equation}
here
$J^\prime(\varphi)(\phi)=\frac{d}{d\lambda}J(\varphi+\lambda\phi)\big|_{\lambda=0}$,
and also $J^\prime(\varphi)(\varphi)=K_1(\varphi)=0$. On the other
hand, a direct computation gives
\begin{align*}K_1^\prime(\varphi)(\varphi)&=2\big(\|\nabla\varphi\|_2^2+\|\varphi\|_2^2\big)-4\iint\limits_{\mathbb{R}^d \mathbb{R}^d}\frac{\varphi(x)^2\varphi(y)^2}{|x-y|^\gamma}dxdy\\
&=-2 \iint\limits_{\mathbb{R}^d
\mathbb{R}^d}\frac{\varphi(x)^2\varphi(y)^2}{|x-y|^\gamma}dxdy<0,\end{align*}
and so $\eta=0$. Hence $J^\prime(\varphi)=0$, namely $
J^\prime(\varphi)(\phi)=0,\ \forall \phi\in H^1:$
$$J^\prime(\varphi)(\phi)=\int_{\mathbb{R}^d}\big(-\Delta\varphi+\varphi-(|x|^{-\gamma}*|\varphi|^2)\varphi\big)\phi dx=0,\
\forall \phi\in H^1.$$ So $\varphi$ satisfies the elliptic equation:
$-\Delta\varphi+\varphi-(|x|^{-\gamma}*|\varphi|^2)\varphi=0$.
Therefore $J(\varphi)\geq J(\mathcal{W})$. On the other hand, it is trivial
that $J(\varphi)\leq J(\mathcal{W})$ by $K_1(\mathcal{W})=0$. Hence $J(\varphi)= J(\mathcal{W})$,
which concludes the proof of Lemma \ref{noempty}.

The following lemma gives an equivalent description of the
functional $K_1(u)$ under  $E(u_0,u_1)<E(\mathcal{W},0).$
\begin{lemma}\label{equi}
Assume $\phi\in H^1$ such that $J(\phi)<J(\mathcal{W})$, then

$(1)\ K_1(\phi)>0\Longleftrightarrow
\|\nabla\phi\|_2^2+\|\phi\|_2^2<\|\nabla \mathcal{W}\|_2^2+\|\mathcal{W}\|_2^2;$

$(2)\ K_1(\phi)<0\Longleftrightarrow
\|\nabla\phi\|_2^2+\|\phi\|_2^2>\|\nabla \mathcal{W}\|_2^2+\|\mathcal{W}\|_2^2.$
\end{lemma}
\begin{proof}
We only prove (1), since (2) can be obtained by the similar way. It
is easy to see that
\begin{equation}\label{J}
J(\mathcal{W})=\frac14\int_{\mathbb{R}^d}\big(|\nabla \mathcal{W}|^2+|\mathcal{W}|^2\big)dx,
\end{equation}
and $$J(\phi)-\frac14 K_1(\phi)=\frac14\int_{\mathbb{R}^d}(|\nabla
\phi|^2+|\phi|^2)dx.$$ Thus, if $J(\phi)<J(\mathcal{W})$ and  $K_1(\phi)>0$,
then
\begin{equation}\label{necc}
\|\nabla\phi\|_2^2+\|\phi\|_2^2<\|\nabla \mathcal{W}\|_2^2+\|\mathcal{W}\|_2^2.
\end{equation}
Next we prove the reverse statement. Since
\begin{equation}\label{suffi}
K_1(e^\lambda\phi)=e^{2\lambda}\Big(\int_{\mathbb{R}^d}[|\nabla\phi|^2+|\phi|^2]dx-e^{2\lambda}
\iint_{\mathbb{R}^d\times
\mathbb{R}^d}\frac{\phi(x)^2\phi(y)^2}{|x-y|^\gamma}dxdy\Big).
\end{equation}
If we take $\lambda\in\mathbb{R}$ such that
\begin{equation}\label{suff1}e^{2\lambda}=\frac{\int_{\mathbb{R}^d}[|\nabla\phi|^2+|\phi|^2]dx}
{\iint_{\mathbb{R}^d\times
\mathbb{R}^d}\frac{\phi(x)^2\phi(y)^2}{|x-y|^\gamma}dxdy},
\end{equation}
then $K_1(e^\lambda\phi)=0.$ By \eqref{suff1}, we know that
$$K_1(\phi)>0\Longleftrightarrow \lambda>0.$$
On the other hand, from $K_1(e^\lambda\phi)=0$ and Lemma
\ref{noempty}, one has \begin{equation}\label{suff2}
J(e^\lambda\phi)\geq J(\mathcal{W}).\end{equation} While by \eqref{suff1} and
the assumption in Lemma \ref{equi}, we obtain
\begin{align}\nonumber
J(e^\lambda\phi)&=e^{2\lambda}\Big\{\frac12\int_{\mathbb{R}^d}\big[|\nabla\phi|^2+|\phi|^2\big]dx-\frac{e^{2\lambda}}{4}
\iint_{\mathbb{R}^d\times \mathbb{R}^d}\frac{\phi(x)^2\phi(y)^2}{|x-y|^\gamma}dxdy\Big\}\\
&=\frac{e^{2\lambda}}{4}\int_{\mathbb{R}^d}\big[|\nabla\phi|^2+|\phi|^2\big]dx\\\nonumber
&<\frac{e^{2\lambda}}{4}\int_{\mathbb{R}^d}\big[|\nabla
\mathcal{W}|^2+|\mathcal{W}|^2\big]dx\\\nonumber&=e^{2\lambda}J(\mathcal{W}).
\end{align}
Combining this with \eqref{suff2}, one gives that $\lambda>0$, and
so $K_1(\phi)>0,$ which concludes the proof.
\end{proof}
\begin{remark}
By the similar argument as above, one may find that if $E(u_0,u_1) <
E(\mathcal{W},0),$ then $K_1(u)\cdot K_2(u) > 0$.
\end{remark}

As a consequence of Lemma \ref{noempty}, we deduce that the sign of
$K_1(u)$ is invariance along the flow of \eqref{equ1} under the
restriction of $E(u,\dot{u})< E(\mathcal{W},0)$:
\begin{corollary}\label{equi1}
Let
\begin{align*}
\mathcal{A}^+&=\{(u,\dot{u})\in H^1\times L^2(\mathbb{R}^d):
E(u,\dot{u})<E(\mathcal{W},0),
K_1(u)>0\},\\
\mathcal{A}^-&=\{(u,\dot{u})\in H^1\times L^2(\mathbb{R}^d):
E(u,\dot{u})<E(\mathcal{W},0), K_1(u)<0\},
\end{align*}
then $\mathcal{A}^+$ and $\mathcal{A}^-$ are invariant under the
flow of \eqref{equ1}. That is, if $(u_0,u_1)\in\mathcal{A}^\pm$,
then $(u(t),\dot{u}(t))\in \mathcal{A}^\pm$, for any
$t\in(-T_-(u_0,u_1),T_+(u_0,u_1))$.
\end{corollary}
\begin{proof}
Since $J(u)<E(u,\dot{u}), J(\mathcal{W})=E(\mathcal{W},0)$ and
$E(u(t),\dot{u}(t))=E(u_0,u_1)$, we deduce that if given
$E(u_0,u_1)<E(\mathcal{W},0)$, then $J(u)<J(\mathcal{W}),$ for any
$t\in(-T_-(u_0,u_1),T_+(u_0,u_1))$.

On the other hand, if there exists $t\in(-T_-(u_0,u_1),T_+(u_0,u_1)$
such that $K_1(u(t))=0$, then by Lemma \ref{noempty}, we get
$J(u(t))\geq J(\mathcal{W})$, which contradicts with $J(u)<J(\mathcal{W})$.

\end{proof}

By Lemma \ref{equi}, one may replace $\mathcal{A}^\pm$ by
\begin{align*}
\mathcal{A}^+&=\{(u,\dot{u})\in H^1\times L^2: E(u,\dot{u})<E(\mathcal{W},0),
\|\nabla u\|_2^2+\|u\|_2^2<\|\nabla \mathcal{W}\|_2^2+\|\mathcal{W}\|_2^2\};\\
\mathcal{A}^-&=\{(u,\dot{u})\in H^1\times L^2: E(u,\dot{u})<E(\mathcal{W},0),
\|\nabla u\|_2^2+\|u\|_2^2>\|\nabla \mathcal{W}\|_2^2+\|\mathcal{W}\|_2^2\}.
\end{align*}

Now, we need a stronger result than Corollary \ref{equi1}, which is
important in the viral analysis. Let
\begin{align*}
\mathcal{A}^+_{\delta,\bar{\delta}}\triangleq\big\{(u,\dot{u})\in
&H^1\times L^2:
E(u,\dot{u})<(1-\delta)E(\mathcal{W},0), K_1(u)>\bar{\delta}\max\{K_1^Q(u), -K_1^N(u)\},\\
& K_2(u)>\bar{\delta}\max\{K_2^Q(u), -K_2^N(u)\}\big\};
\end{align*}
\begin{proposition}\label{conservation}
Assume $(\phi_0,\phi_1)\in H^1\times L^2$ and $(\phi_0,\phi_1)\in
\mathcal{A}^+_{\delta,0}$, then there exists some
$\bar{\delta}\in(0,1)$ independent of $(\phi_0,\phi_1)$ such that
$(\phi_0,\phi_1)\in \mathcal{A}^+_{\delta,\bar{\delta}}.$
\end{proposition}

\begin{proof}
The proof is the same as \cite{MW}. For convenience, we give a full
proof. We only need to prove
$$K_1(\phi_0)>\bar{\delta}\max\big\{\|\nabla
\phi_0\|_2^2+\|\phi_0\|_2^2, -K_1^N(\phi_0)\big\},$$ since the other
is given by the similar way. From $(\phi_0,\phi_1)\in
\mathcal{A}^+_{\delta,0}$, we know that $K_1(\phi_0)>0$, so it
suffices to prove that
\begin{equation}\label{su1}
K_1(\phi_0)>\bar{\delta}\big(\|\nabla
\phi_0\|_2^2+\|\phi_0\|_2^2\big).
\end{equation}
We suppose for contradiction that there exists a sequence
$\{(\varphi_0^n,\varphi_1^n)\in\mathcal{A}^+_{\delta,0}\}_n$, and a
sequence $\{\bar{\delta}_n\}_n$ satisfying
$\bar{\delta}_n\rightarrow 0,$ as $n\rightarrow\infty$, and
$$J(\varphi_0^n)<(1-\delta)J(\mathcal{W}), 0<K_1(\varphi_0^n)\leq\bar{\delta}_n(\|\nabla
\varphi_0^n\|_2^2+\|\varphi_0^n\|_2^2),
$$
which implies that
$$|K_1^N(\varphi_0^n)|<K_1^Q(\varphi_0^n)\leq\frac{1}{1-\bar{\delta}_n}|K_1^N(\varphi_0^n)|.$$
Thus there exists $\lambda_n$ such that
$$K_1(e^{\lambda_n}\varphi_0^n)=0,$$ and
$1<e^{2\lambda_n}\leq\frac{1}{1-\bar{\delta}_n}.$ So by Lemma
\ref{noempty}, we get $J(e^{\lambda_n}\varphi)\geq J(\mathcal{W})$, which
implies that \begin{equation}\label{min1}\frac{e^{2\lambda_n}}{2}
K_1^Q(\varphi_0^n)+\frac{e^{4\lambda_n}}{4} K_1^N(\varphi_0^n)\geq
J(\mathcal{W}).\end{equation} On the other hand, since
$(\varphi_0^n,\varphi_1^n)\in\mathcal{A}^+_{\delta,0}$,
$J(\varphi_0^n)<E(\varphi_0^n,\varphi_1^n)<(1-\delta)E(\mathcal{W},0)=(1-\delta)J(\mathcal{W})$,
which means that \begin{equation}\label{min2}\frac{1}{2}
K_1^Q(\varphi_0^n)+\frac{1}{4} K_1^N(\varphi_0^n)<
(1-\delta)J(\mathcal{W}).\end{equation} Combining \eqref{min1} and
\eqref{min2}, one gives
$$0\leq-\frac{e^{2\lambda_n}-1}{4}K_1^N(\varphi_0^n)<(1-\delta-e^{-2\lambda_n}) J(\mathcal{W}).$$
But this can not happen for large $n$, since
$e^{2\lambda_n}\rightarrow 1$ as $n\rightarrow\infty$.
\end{proof}

Combining the energy conservation law, Corollary \ref{equi1} with
the above proposition, we obtain the following result.

\begin{corollary}\label{reserve}
Assume $(u_0,u_1)\in\mathcal{A}^+$, then there exist some
$\delta,\bar{\delta}\in (0,1)$ depending on $(u_0,u_1)$, such that
the corresponding solution
$(u(t),\dot{u}(t))\in\mathcal{A}^+_{\delta,\bar{\delta}}$ for any
$t\in(-T_-(u_0,u_1),T_+(u_0,u_1)).$
\end{corollary}

It is easy to observe that the free energy and the nonlinear energy
are equivalent in the set $K_1> 0$.
\begin{lemma}\label{freeequi}
Assume $2<\gamma<\min\{4,d\},$ then for any $(u_0,u_1)\in H^1\times
L^2$, we have
\begin{equation}
K_1(u_0)> 0\ \Longrightarrow\ \begin{cases} J(u_0)\leq\frac12
K_1^Q(u_0)\leq 2J(u_0)\\
E(u_0,u_1)\leq E_0(u_0,u_1)\leq2 E(u_0,u_1)\end{cases}
\end{equation}
where $E_0(u_0,u_1)=\frac12 K_1^Q(u_0)+\frac12\|u_1\|_2^2$, and
$E(u_0,u_1)=E_0(u_0,u_1)+\frac14 K_1^N(u_0).$
\end{lemma}
\begin{proof}
First, we recall that
\begin{equation*}
\begin{cases}
K_1(u_0)=K_1^Q(u_0)+K_1^N(u_0),\\
J(u_0)=\frac12K_1^Q(u_0)+\frac14K_1^N(u_0).
\end{cases}
\end{equation*}
Hence, by direct computation, we get $K_1(u_0)> 0,$ which gives
$J(u_0)\leq\frac12 K_1^Q(u_0)\leq 2J(u_0).$ Similarly, we obtain
$E(u_0,u_1)\leq E_0(u_0,u_1)\leq2 E(u_0,u_1)$ under the condition
$K_1(u_0)> 0$.
\end{proof}
\begin{corollary}\label{global}
Let $\mu=-1,~2<\gamma<\min\{d,4\}$ and $(u_0,u_1)\in
H^1(\mathbb{R}^d\times L^2(\mathbb{R}^d)$. Assume that
$E(u_0,u_1)<E(\mathcal{W},0)$ and $K_1(u_0)>0$, then the solution of
\eqref{equ1} is global.
\end{corollary}
\begin{proof}
It is a consequence of Corollary \ref{equi1}, energy conservation,
Lemma \ref{freeequi} and the standard blow-up criterion, see Theorem
\ref{small}.
\end{proof}

The next lemma gives a upper bound on $K_1$ in the set $K_1<0$,
which will be important for the blow up.
\begin{lemma}\label{convex}
Suppose $2<\gamma<\min\{4,d\}$, $\phi\in H^1$, $J(\phi)<J(\mathcal{W})$, and
$K_1(\phi)<0$, then
$$K_1(\phi)<-2(J(\mathcal{W})-J(\phi)).$$

\end{lemma}
\begin{proof}
Let $j(\lambda)=J(e^\lambda\phi),$ then by a directive computation,
we have
$$j^\prime(\lambda)=K_1(e^\lambda\phi),\
 j^{\prime\prime}(\lambda)=2K_1^Q(e^\lambda\phi)+4K_1^N(e^\lambda\phi),$$
and so $j^\prime(0)=K_1(\phi)$,
$j^{\prime\prime}(\lambda)<2j^\prime(\lambda)$. If we choose
$\lambda_0\in\mathbb{R}$ such that
$$e^{2\lambda_0}=\frac{K_1^Q(\phi)}{-K_1^N(\phi)},$$
then $j^\prime(\lambda_0)=K_1(e^{\lambda_0}\phi)=0$, and
$\lambda_0<0$ by $K_1(\phi)<0$. Therefore, we have
\begin{equation}\label{connn}
K_1(\phi)=j^\prime(0)-j^\prime(\lambda_0)=\int_{\lambda_0}^0j^{\prime\prime}
(\lambda)d\lambda<2\int_{\lambda_0}^0j^{\prime}(\lambda)d\lambda=2(j(0)-j(\lambda_0)).
\end{equation}
By Lemma \ref{noempty} and $K_1(e^{\lambda_0}\phi)=0$, we get
$j(\lambda_0)=J(e^{\lambda_0}\phi)\geq J(\mathcal{W})$, and so
$$\eqref{connn}\leq 2(J(\phi)-J(\mathcal{W})).$$
This completes the proof.
\end{proof}
\begin{corollary}\label{blow-up1}
Suppose $2<\gamma<\min\{4,d\}$, $(u_0,u_1)\in H^1\times L^2$,
$E(u_0,u_1)<E(\mathcal{W},0)$ and $K_1(u_0)<0$, then
$$K_1(u)\leq -2(E(\mathcal{W},0)-E(u_0,u_1)),~\ \text{for\ any}\  ~ t\in
(-T_-(u_0,u_1),T_+(u_0,u_1)).$$
\end{corollary}
\begin{proof}
From Corollary \ref{equi1}, we have $$K_1(u)<0, \ \text{for}\ t\in
(-T_-(u_0,u_1),T_+(u_0,u_1)).$$ And so by Lemma \ref{convex}, we
obtain $$K_1(u)\leq -2(E(\mathcal{W},0)-E(u_0,u_1)),\ \text{for}\ t\in
(-T_-(u_0,u_1),T_+(u_0,u_1)),$$ where we use the fact that $J(u)\leq
E(u,\dot{u})=E(u_0,u_1),$ and $J(\mathcal{W})=E(\mathcal{W},0)$.
\end{proof}

\section{Blow up}
\setcounter{section}{4}\setcounter{equation}{0} In this section we
prove the blow-up part of Theorem \ref{theorem1.2}. The idea is
essentially due to Payne-Sattinger \cite{PS}, but we give a complete
proof for convenience.

By contradiction we assume that the solution $u$ exists for all
$t>0$. The proof for $t<0$ is the same.

 Denote $y(t)=\int|u|^2dx$,
then we have $y'(t)=2\int uu_tdx,$ and
\begin{equation}\label{blowup}y''(t)=2\|\dot{u}\|_2^2-2K_1(u)=6\|\dot{u}\|_2^2+2K_1^Q(u)-8E(u_0,u_1).\end{equation}
It follows from Corollary \ref{blow-up1} that there exists
$\delta>0$ such that $K_1(u)<-\delta$, and so $y''(t)>2\delta$. Thus
$$y''(t)>2\delta>0.$$
 Then by the lower bound on $y''(t)$, there
exists $t_0>0$ such that $y'(t_0)>0$, and hence $y'(t)>0$ for
$t>t_0$. By Lemma \ref{equi} and $K_1(u)<0$, we get
$K_1^Q(u)>K_1^Q(\mathcal{W}),$ and so
$$E(u,\dot{u})<E(\mathcal{W},0)=J(\mathcal{W})=\frac14K_1^Q(\mathcal{W})<\frac14K_1^Q(u).$$ Therefore, using Cauchy-Schwarz inequality and \eqref{blowup}, we obtain for
any $t>t_0$
$$y''(t)\geq 6\|\dot{u}\|_2^2\geq\frac32\frac{y'^2}{y}.$$
So that, for $t>t_0$,
$$\frac{y''(t)}{y'(t)}\geq\frac32\frac{y'(t)}{y(t)}.$$
Hence for $t>t_0$,
$$y(t)^{-\frac12}\leq
y(t_0)^{-\frac12}-\frac12\frac{y'(t_0)}{y(t_0)^\frac32}(t-t_0).$$
Therefore $T_+\leq t_0+\frac{2y'(t_0)}{y(t_0)^{\frac32}},$ which
contradicts with $T_+=+\infty$.

\section{Profile decomposition}
\setcounter{section}{5}\setcounter{equation}{0}

In this section, we first recall the linear profile decomposition of
the sequence of $H^1$-bounded solutions of \eqref{equ1} which was
established in \cite{IMN}. And then we utilize it to show the
orthogonal analysis for the nonlinear energy and the nonlinear
profile decomposition which will be used to construct the critical
element and obtain its compactness properties. In order to do it, we
now recall some notations in \cite{IMN}.

 With any real-valued function
$u(t,x)$, we associate the complex-valued function $\vec{u}(t,x)$ by
\begin{equation}
\vec{u}=\langle\nabla\rangle u-i\dot{u},\quad
u=\Re\langle\nabla\rangle^{-1}\vec{u}.
\end{equation}
Then the free and nonlinear Klein-Gordon equations are given by
\begin{align}
\begin{cases}
(\Box+1)u=0\Longleftrightarrow(i\partial_t+\langle\nabla\rangle)\vec{u}=0,\\
(\Box+1)u=f(u)\Longleftrightarrow(i\partial_t+\langle\nabla\rangle)\vec{u}=f(\langle\nabla\rangle^{-1}\Re\vec{u}),
\end{cases}
\end{align}
and the energy are written as
$$\tilde{E}(\vec{u})=E(u,\dot{u})=\frac12 \int_{\mathbb{R}^d} \big(\big|\dot{u}
\big|^2 +\big| \nabla u \big|^2 + \big| u \big|^2 \big)dx +
\frac{\mu}{4} \iint_{\mathbb{R}^d\times\mathbb{R}^d}
 \frac{|
u(t,x)|^2 |u(t,y)|^2}{|x-y|^\gamma} dxdy.$$
 We denote the set of Fourier multipliers on
$$\mathcal{MC}=\{\mu=\mathcal{F}^{-1}\tilde{\mu}\mathcal{F}|\ \tilde{\mu}\in
C(\mathbb{R}^d),\exists
\lim\limits_{|x|\rightarrow\infty}\tilde{\mu}(x)\in\mathbb{R}\}.$$

 \subsection{Linear profile decomposition}
First, we state the linear profile decomposition(which was
established in\cite{IMN}) as follows
\begin{lemma}\label{lem3.1}
Let $\vec{v}_n$ be a sequence of free Klein-Gordon solutions with
uniformly bounded $L^2_x$ norm. Then after replacing it with some
subsequence, there exist $K\in\{0,1,2\ldots,\infty\}$ and, for each
integer $j\in[0,K)$, $\varphi^j\in L^2(\mathbb{R}^d)$ and $\{(t^j_n,
x^j_n)\}_{n\in\mathbb{N}}\subset\mathbb{R}\times\mathbb{R}^d$
satisfying the following. Define $\vec{v}^j_n$ and
$\vec{\omega}^k_n$ for each $j<k\leq K$ by
\begin{equation}\label{equ3.1}
\vec{v}_n(t,x)=\sum\limits_{j=0}^{k-1}\vec{v}^j_n(t,x)+\vec{\omega}^k_n(t,x),\quad
\vec{v}^j_n(t,x)=e^{i\langle\nabla\rangle(t-t^j_n)}\varphi^j(x-x^j_n),
\end{equation}
then for any $s<-\frac{d}{2}$, we have
\begin{equation}\label{equ3.2}
\lim\limits_{k\rightarrow
K}\varlimsup\limits_{n\rightarrow\infty}\|\vec{\omega}^k_n\|_{L^\infty(\mathbb{R};B^s_{\infty,1}(
\mathbb{R}^d))}=0,
\end{equation}
and for any $\mu\in \mathcal{MC}$, any $l<j<k\leq K$ and any $t\in
\mathbb{R}$,
\begin{align}\label{equ3.3}
\lim\limits_{n\rightarrow\infty}\langle\mu\vec{v}^l_n,
\mu\vec{v}^j_n\rangle_{L^2_x}^2=0=
\lim\limits_{n\rightarrow\infty}\langle\mu\vec{v}^j_n,
\mu\vec{\omega}^k_n\rangle_{L^2_x}^2,
\\\label{equ3.4}
\lim\limits_{n\rightarrow\infty}|t_n^j-t_n^k|+|x_n^j-x_n^k|=+\infty.
\end{align}
\end{lemma}
\begin{remark}
We call $\vec{v}_n^j$ the free concentrating wave. From
\eqref{equ3.3}, we have the following asymptotic orthogonality
\begin{equation}
\lim\limits_{n\rightarrow+\infty}\Big(\|\mu\vec{v}_n(t)\|_{L^2}^2-\sum\limits_{j=0}^{k-1}\|\mu\vec{v}_n^j(t)\|_{L^2}^2
-\|\mu\vec{\omega}_n^k(t)\|_{L^2}^2\Big)=0.
\end{equation}
\end{remark}
Next we begin with the orthogonal analysis for the nonlinear energy.
\begin{lemma}\label{energy}
Let $\vec{v}_n$ be a sequence of free Klein-Gordon solutions
satisfying $\vec{v}_n(0)\in L^2_x$. Let
$\vec{v}_n=\sum\limits_{j=0}^{k-1}\vec{v}^j_n+\vec{\omega}_n^k$ be
the linear profile decomposition given by Lemma \ref{lem3.1}. Then
if $\mu=1$ (defocusing)  and
$\varlimsup\limits_{n\rightarrow\infty}\tilde{E}(\vec{v}_n(0))<+\infty$,
then we have $\vec{v}_n^j(0)\in L^2_x$ for large $n$, and
\begin{equation}\label{equ3.5}
\lim\limits_{k\rightarrow
K}\varlimsup\limits_{n\rightarrow\infty}\Big|\tilde{E}(\vec{v}_n(0))-\sum\limits_{j=0}^{k-1}
\tilde{E}(\vec{v}^j_n(0))-\tilde{E}(\vec{\omega}^k_n(0))\Big|=0.
\end{equation}
Moreover we have for all $j<k$
\begin{equation}\label{equ3.6}
0\leq\varliminf\limits_{n\rightarrow\infty}\tilde{E}(\vec{v}^j_n(0))\leq\varlimsup\limits_{n\rightarrow
\infty}\tilde{E}(\vec{v}^j_n(0))\leq\varlimsup\limits_{n\rightarrow
\infty}\tilde{E}(\vec{v}_n(0)),
\end{equation}
where the last inequality becomes equality only if $K=1$ and
$\vec{\omega}^1_n\rightarrow 0$ in $L^\infty_t L^2_x$.

Furthermore, if $\mu=-1$ (focusing), $(v_n(0),\dot{v}_n(0))\in
\mathcal{A}^+$ and
$\varlimsup\limits_{n\rightarrow\infty}\tilde{E}(\vec{v}_n(0))<E(\mathcal{W},0)$,
then we have $\vec{v}_n^j(0)\in L^2_x$ for large $n$ and all $j<K$,
$(v_n^j(0),\dot{v}_n^j(0))\in\mathcal{A}^+$, and \eqref{equ3.5},
\eqref{equ3.6} also holds true.
\end{lemma}

\begin{proof}
By Sobolev imbedding theorem and \eqref{equ3.2}, we have
$$\lim\limits_{k\rightarrow
K}\varlimsup\limits_{n\rightarrow\infty}\|\omega^k_n\|_{L_x^\frac{4d}{2d-\gamma}}\leq\lim\limits_{k\rightarrow
K}\varlimsup\limits_{n\rightarrow\infty}\|\omega^k_n\|_{\dot{H}^1}^\theta\|\omega^k_n\|_{B^{\frac{1-d}{2}}_{\infty,1}}^{1-\theta}=0,$$
where $\theta=\frac{\gamma-2}{2}\in (0,1)$, and $\omega^k_n=\Re
\langle\nabla\rangle^{-1}\vec{\omega}^k_n$. This implies that, if
there exists $u_i=\omega_n^k\ (i=1,2,3,4)$, then by the H\"{o}lder
and general Young inequality, we obtain
\begin{align*}
\lim\limits_{k\rightarrow
K}\varlimsup\limits_{n\rightarrow\infty}\|\big(|x|^{-\gamma}*(u_1u_2)\big)(u_3u_4)\|_{L^1_x}\leq
\lim\limits_{k\rightarrow
K}\varlimsup\limits_{n\rightarrow\infty}\prod\limits_{i=1}^{4}\|u_i\|_{L_x^\frac{4d}{2d-\gamma}}
=0.
\end{align*}
This together with \eqref{equ3.3} reduces us to prove
$$\lim\limits_{k\rightarrow
K}\varlimsup\limits_{n\rightarrow\infty}\Big\{\big\|(|x|^{-\gamma}*|v_n-\omega_n^k|^2)|v_n-\omega_n^k|^2\big\|_{L_x^1}-\sum\limits_{j=0}^{k-1}\big\|(|x|^{-\gamma}*|v_n^j|^2)|v_n^j|^2\big\|_{L_x^1}\Big\}=0.$$
For this purpose, we discuss in two  cases
\begin{align*}
\begin{cases}
\text{case\ 1:}\ \exists\ j,\ |t_n^j|\rightarrow\infty\ \text{as}\
n\rightarrow\infty,\\
\text{case\ 2:}\ |t_n^j|\ \text{is\ uniformly\ bounded\ in}\ n, j.
\end{cases}
\end{align*}

For the first case, by the decay of $e^{it\langle\nabla\rangle}$ in
$\mathcal{S}\rightarrow L^p_x$ uniform w.r.t $n$ and the Sobolev
embedding $\dot{H}^s\subset L^p$, we have
$$\|v_n^j\|_{L_x^\frac{4d}{2d-\gamma}}=\|\Re\ \langle\nabla\rangle^{-1}e^{i\langle\nabla\rangle(t-t^j_n)}\varphi^j(x-x^j_n)\|_{L_x^\frac{4d}{2d-\gamma}}\rightarrow
0,\ \text{as}\ n\rightarrow\infty.$$  Thus by the linear profile
decomposition, the H\"{o}lder and generalized Young inequality, we
have
\begin{align*}
\big\|\big(|x|^{-\gamma}*(v_n^{k_1}v_n^{k_2})\big)(v_n^{k_3}v_n^{k_4})\big\|_{L_x^1}\lesssim\prod\limits_{i=1}^{4}
\big\|v_n^{k_i}\big\|_{L_x^\frac{4d}{2d-\gamma}}\rightarrow 0, \
\text{as}\ n\rightarrow\infty,
\end{align*}
for some $k_i=j$.

Next we consider the second case. Since
\begin{align*}
\Big\|(|x|^{-\gamma}*|\sum\limits_{j=0}^{k-1}v_n^j|^2)|\sum\limits_{j=0}^{k-1}v_n^j|^2\Big\|_{L_x^1}=\sum\limits_{j_1,j_2,j_3,j_4<k}\Big\|
(|x|^{-\gamma}*(v_n^{j_1}v_n^{j_2}))|v_n^{j_3}v_n^{j_4}|\Big\|_{L_x^1},
\end{align*}
we only need to prove that
\begin{equation}\label{equ31}\Big\|(|x|^{-\gamma}*(v_n^{j_1}v_n^{j_2}))|v_n^{j_3}v_n^{j_4}|\Big\|_{L_x^1}\rightarrow
0,\ \text{as}\ n\rightarrow\infty,\end{equation} provided that
$t_n^{j_i}$ is bounded for any $i=1,2,3,4$, and at least two of
$j_1,j_2,j_3,j_4$ are different. Moreover, by \eqref{equ3.4}, we
know that
\begin{equation}\label{orthong}|x_n^{j_i}-x_n^{j_l}|\rightarrow\infty,\
\text{as}\ n\rightarrow\infty,\ \text{for}\
j_i\not=j_l.\end{equation}

To prove \eqref{equ31}, we should split it into the following two
cases:
\begin{align*}
\begin{cases}
\text{Case\ 1:}\ j_1\neq j_2\ \text{or} \ j_3\neq j_4,\\
\text{Case\ 2:}\ j_1=j_2,\ j_3=j_4,\ j_1\neq j_3.
\end{cases}
\end{align*}

Keep in mind that
$\vec{v}^j_n(t,x)=e^{i\langle\nabla\rangle(t-t^j_n)}\varphi^j(x-x^j_n)$.
Without loss of generality, we may assume that
$\varphi^{j_i}(i=1,2,3,4)$ have compact support in $\{x\in
\mathbb{R}^d: |x|\leq R\}$.

\noindent{\bf Case1: $j_1\neq j_2$ or $j_3\neq j_4$,}
$$
\text{LHS}\ \text{of}\ \eqref{equ31}\lesssim
\|v_n^{j_1}v_n^{j_2}\|_p\|v_n^{j_3}v_n^{j_4}\|_p\rightarrow 0,\quad
\text{as}\ n\rightarrow\infty,
$$
where we use the orthogonal condition \eqref{orthong} and
$p=\frac{2d}{2d-\gamma}$.

\noindent{\bf Case2: $j_1=j_2$, $j_3=j_4$, $j_1\neq j_3$.}
\begin{align*}
\text{LHS}\ \text{of}\
\eqref{equ31}&=\|(|x|^{-\gamma}*|v_n^{j_1}|^2)
|v_n^{j_3}|^2\|_{L_x^1}\\
&\leq\|(V_1*|v_n^{j_1}|^2)|v_n^{j_3}|^2\|_{L_x^1}+\|(V_2*|v_n^{j_1}|^2)|v_n^{j_3}|^2\|_{L_x^1}\\
&\triangleq I_1+I_2.
\end{align*}
where $V_1=|x|^{-\gamma}\chi_{|x|\leq LR},
V_2=|x|^{-\gamma}\chi_{|x|\geq LR}$ for large $L\gg 1$. By supp
$\varphi^{j_1}\subset\{x\in \mathbb{R}^d: |x|\leq R\}$, the H\"older
inequality and Young inequality, we get
\begin{align}\nonumber
I_1&=\big\|(V_1*||e^{-it^{j_1}_n\langle\nabla\rangle}\varphi^{j_1}(\cdot-x^{j_1}_n)|^2)
|e^{-it^{j_3}_n\langle\nabla\rangle}\varphi^{j_3}(x-x^{j_3}_n)|^2\big\|_{L_x^1}\\\nonumber
&=\Big\|\big[V_1*|e^{-it^{j_1}_n\langle\nabla\rangle}\varphi^{j_1}(\cdot)|^2\big]\chi_{|x|\lesssim
R}(x)|e^{-it^{j_3}_n\langle\nabla\rangle}\varphi^{j_3}\big(x-(x_n^{j_3}-x_n^{j_1})\big)|^2\Big\|_{L_x^1}\\\nonumber
&\leq
\big\|V_1*|e^{-it^{j_1}_n\langle\nabla\rangle}\varphi^{j_1}(\cdot)|^2\big\|_{L^\frac{2d}{\gamma}}\big\|\chi_{|x|\lesssim
R}(x)e^{-it^{j_3}_n\langle\nabla\rangle}\varphi^{j_3}\big(x-(x_n^{j_3}-x_n^{j_1})\big)\big\|_{L_x^\frac{4d}{2d-\gamma}}^2\\\nonumber
&\lesssim
\|V_1\|_{L_*^{\frac{d}{\gamma}}}\|v_{n}^{j_1}\|_{L^\frac{4d}{2d-\gamma}}^2\big\|\chi_{|x|\lesssim
R}(x)e^{-it^{j_3}_n\langle\nabla\rangle}\varphi^{j_3}\big(x-(x_n^{j_3}-x_n^{j_1})\big)\big\|_{L_x^\frac{4d}{2d-\gamma}}^2\\\label{equat1}
&\rightarrow 0, \quad \text{as}\quad n\rightarrow \infty.
\end{align}
 For $I_2$, by the H\"older and Young
inequality, one infers that
\begin{align}\nonumber
I_2&=\big\|(V_2*||v_{n}^{j_1}|^2)
|v_n^{j_3}|^2\big\|_{L_x^1}\\\nonumber
&\leq\|V_2*|v_{n}^{j_1}|^2\|_{L^\infty}\|v_n^{j_3}\|_{L^2}^2\\\nonumber
&\lesssim
(LR)^{-\gamma}\|v_{n}^{j_1}\|_{L^2}^2\|v_n^{j_3}\|_{L^2}^2\\\label{equat2}
&\rightarrow 0, \quad \text{as} \quad L\rightarrow\infty.
\end{align}
Combining \eqref{equat1} with \eqref{equat2}, we obtain
\eqref{equ31},  which concludes the proof of \eqref{equ3.5}.

Now we turn to prove that
$(v_n^j(0),\dot{v}_n^j(0))\in\mathcal{A}^+$ for large $n$ and all
$j<K$ in the focusing case. In fact, by \eqref{equ3.3}, we have
$$\lim\limits_{k\rightarrow
K}\varlimsup\limits_{n\rightarrow\infty}\bigg|K_1^Q(v_n(0))-\sum\limits_{j<k}K_1^Q(v_n^j(0))-K_1^Q(\omega_n^k(0))\bigg|=0,$$
where $K_1^Q(\phi)=\|\phi\|_2^2+\|\nabla\phi\|_2^2.$ It follows from
 $v_n(0)\in \mathcal{A}^+$ and Lemma \ref{equi} that
$$K_1^Q(v_n^j(0))<K_1^Q(\mathcal{W}),$$
and so $(v_n^j(0),\dot{v}_n^j(0))\in\mathcal{A}^+$ for large $n$.
Lemma \ref{freeequi} shows that the last
inequality in \eqref{equ3.6} becomes equality only if $K=1$ and
$\vec{\omega}^1_n\rightarrow 0$ in $L^\infty_t L^2_x$.
\end{proof}

\subsection{Nonlinear profile decomposition}
After the linear profile decomposition of a sequence of initial data
in the last subsection, we now show the nonlinear profile
decomposition of a sequence of the solutions of \eqref{equ1} with
the same initial data in the energy space $H^1(\mathbb{R}^d)\times
L^2(\mathbb{R}^d).$

First we construct a nonlinear profile corresponding to a free
concentrating wave. Let $\vec{v}_n$ be a free concentrating wave for
a sequence $(t_n,x_n)$,
\begin{align}
\begin{cases}
(i\partial_t+\langle\nabla\rangle)\vec{v}_n=0,\\
\vec{v}_n(t_n)=\phi(x-x_n),\ \phi(x)\in L^2,
\end{cases}
\end{align}
and let $u_n$ be the nonlinear solution with the same initial data
\begin{align}
\begin{cases}
(i\partial_t+\langle\nabla\rangle)\vec{u}_n=f(u_n),\\
\vec{u}_n(0)=\vec{v}_n(0).
\end{cases}
\end{align}
Next we define
$$\vec{V}_n(t,x)=\vec{v}_n(t+t_n,x+x_n),~\vec{U}_n(t,x)=\vec{u}_n(t+t_n,x+x_n).$$
Then they satisfy the rescaled equations
$$\vec{V}_n(t,x)=e^{it\langle\nabla\rangle}\phi(x),\quad\vec{U}_n(t)=\vec{V}_n(t)-i\int_{-t_n}^{t}e^{i(t-s)\langle\nabla\rangle}f(\Re\langle\nabla\rangle^{-1}\vec{U}_n)ds.$$
Extracting a subsequence, we may assume convergence
$$t_n\rightarrow t_\infty\in[-\infty,\infty].$$
Thus the limit equations are given by
$$\vec{V}_\infty(t,x)=e^{it\langle\nabla\rangle}\phi(x),\quad\vec{U}_\infty(t)=\vec{V}_\infty(t)-i\int_{-t_\infty}^{t}e^{i(t-s)\langle\nabla\rangle}f(\Re\langle\nabla\rangle^{-1}\vec{U}_\infty)ds.$$
The unique existence of a local solution $\vec{U}_\infty$ around
$t=t_\infty$ is know in all cases, including $t=\pm\infty$ which
corresponds to the existence of the wave operators, by using the
standard iteration with the Strichartz estimate.

\begin{definition}\label{def3.1}
The nonlinear concentrating wave $\vec{u}_{(n)}$ associated with
$\vec{v}_n$ is defined by
$$\vec{u}_{(n)}(t,x)=\vec{U}_\infty(t-t_n,x-x_n).$$

\end{definition}

\begin{remark}\label{remark5.2}
(i) $u_{(n)}$ solves \eqref{equ1}.

(ii) By definition, we have
$\|\vec{u}_n(0)-\vec{u}_{(n)}(0)\|_{L^2_x}\rightarrow 0,$ as
$n\rightarrow \infty.$ In fact
\begin{align*}\nonumber
\|\vec{u}_n(0)-\vec{u}_{(n)}(0)\|_{L^2_x}&=\|\vec{v}_n(0)-\vec{u}_{(n)}(0)\|_{L^2_x}\\\nonumber
&=\|\vec{V}_n(-t_n)-\vec{U}_\infty(-t_n)\|_{L^2_x}\\\nonumber &=
\Big\|\vec{U}_\infty(-t_\infty)-\vec{U}_\infty(-t_n)\|\Big\|_{L^2_x}\\
&\rightarrow 0,\quad \text{as}\quad n\rightarrow\infty.
\end{align*}
\end{remark}
Let $u_n$ be a sequence of (local) solutions of \eqref{equ1} around
$t=0$, and let $v_n$ be the sequence of the free solutions with the
same initial data. We consider the linear profile decomposition of
$\{\vec{v}_n\}$ given by Lemma \ref{lem3.1},
\begin{equation*}
\vec{v}_n=\sum\limits_{j=0}^{k-1}\vec{v}^j_n+\vec{\omega}^k_n,\quad
\vec{v}^j_n=e^{i\langle\nabla\rangle(t-t^j_n)}\varphi^j(x-x^j_n).
\end{equation*}

\begin{definition}
(Nonlinear profile decomposition) Let
$\{\vec{v}_n^j\}_{n\in\mathbb{N}}$ be the free concentrating wave,
and $\{\vec{u}_{(n)}^j\}_{n\in\mathbb{N}}$ be the sequence of the
nonlinear concentrating wave associated with
$\{\vec{v}_n^j\}_{n\in\mathbb{N}}$. Then we define the nonlinear
profile decomposition of $u_n$ by
\begin{equation}\label{nonlineard}
\vec{u}_{(n)}^{<k}:=\sum\limits_{j=0}^{k-1}\vec{u}_{(n)}^j.
\end{equation}
\end{definition}

We are going to prove that $\vec{u}_{(n)}^{<k}+\vec{\omega}_n^k$ is
a good approximation for $\vec{u}_n$. And the following two lemmas derive
from Lemma \ref{lem3.1} and the perturbation lemma. The first lemma
concerns the orthogonality in the Strichartz norms.
\begin{lemma}\label{lem3.3}
Suppose that in the nonlinear profile decomposition
\eqref{nonlineard}, we have
\begin{equation}\label{equ3.7}
\|\Re\langle\nabla\rangle^{-1}\vec{U}^j_\infty\|_{ST(\mathbb{R})}+\|\vec{U}^j_\infty\|_{L_t^\infty
L^2_x(\mathbb{R})}<\infty,\ \forall j<k.
\end{equation}
Then, for any finite interval $I, j<k,$ we have
\begin{align}\label{equ3.8}
\varlimsup\limits_{n\rightarrow\infty}\|u_{(n)}^j\|_{ST(I)}&\lesssim\|\Re\langle\nabla\rangle^{-1}\vec{U}_\infty^j\|_{ST(\mathbb{R})},\\\label{equ3.9}
\varlimsup\limits_{n\rightarrow\infty}\|u_{(n)}^{<k}\|_{ST(I)}^2&\lesssim\varlimsup\limits_{n\rightarrow\infty}
\sum\limits_{j=0}^{k-1}\big\|u_{(n)}^j\big\|_{ST(\mathbb{R})}^2,
\end{align}
where the implicit constants do not depend on $I$ or $j$. We also
have
\begin{equation}\label{equ3.10}
\lim\limits_{n\rightarrow\infty}\bigg\|(|x|^{-\gamma}*|u_{(n)}^{<k}|^2)u_{(n)}^{<k}-\sum\limits_{j=0}^{k-1}(|x|^{-\gamma}*|u_{(n)}^j|^2)u_{(n)}^j\bigg\|_{ST^*(I)}=0.
\end{equation}
\end{lemma}

\begin{proof}
It is easy to get \eqref{equ3.8} from the definition of $u_{(n)}^j$.
Now we prove \eqref{equ3.9}. Let $\chi(t,x)\in
C^\infty(\mathbb{R}^{n+1})$ satisfy $\chi(t,x)=1$ for $|(t,x)|\leq
1$ and $\chi(t,x)=0$ for $|(t,x)|\geq 2$, and
$$\chi_R(t,x)=\chi(\frac{t}{R},\frac{x}{R}).$$
Define  $u_{(n),R}^j$ for $R\gg1$ by
$$
u_{(n),R}^j(t,x)=\chi_R(t,x)u_{(n)}^j(t,x),\quad
u_{(n),R}^{<k}=\sum\limits_{j=0}^{k-1}u_{(n),R}^j.$$ Then we have
$$\|u_{(n)}^{<k}-u_{(n),R}^{<k}\|_{ST(I)}\leq
\sum\limits_{j=0}^{k-1}\big\|(1-\chi_R)\Re\langle\nabla\rangle^{-1}\vec{U}_\infty^j\big\|_{ST(\mathbb{R})}\rightarrow
0,\quad \text{as}\ R\rightarrow\infty,$$ thus we may replace
$u_{(n)}^{<k}$ by $u_{(n),R}^{<k}$. The homogenous Besov norm
in ST is equivalent to
\begin{align}\label{equi2}
\|u_{(n),R}^{<k}\|_{L^q(I;B^{\frac12}_{q,2})}&\cong\|u_{(n),R}^{<k}\|_{L^q(I;\dot{B}^{\frac12}_{q,2})}+\|u_{(n),R}^{<k}\|_{L^q_{t,x}}\\\nonumber
&\cong\bigg(\int_{I}\Big(\int_{\mathbb{R}^d}\frac{\|u_{(n),R}^{<k}(x-y)-u_{(n),R}^{<k}(x)\|_{L^q}^2}{|y|}\frac{dy}{|y|^d}\Big)
^{\frac{q}{2}}\bigg)^{\frac1q}+\|u_{(n),R}^{<k}\|_{L^q_{t,x}},
\end{align}
By the orthogonality \eqref{equ3.4}, we  get for large $n$
$$|u_{(n),R}^{<k}(x-y)-u_{(n),R}^{<k}(x)|=\bigg\{\sum\limits_{j=0}^{k-1}[u_{(n),R}^{j}(x-y)-u_{(n),R}^{j}(x)]^2\bigg\}^\frac12,$$
and so
\begin{align*}
\eqref{equi2}&=\bigg(\int_{I}\Big(\int_{\mathbb{R}^d}\frac{\Big\|\Big\{\sum\limits_{j=0}^{k-1}[u_{(n),R}^{j}(x-y)-u_{(n),R}^{j}(x)]^2
\Big\}^\frac12\Big\|_{L^q}^2}{|y|}\frac{dy}{|y|^d}\Big)
^{\frac{q}{2}}\bigg)^{\frac1q}+\|u_{(n),R}^{<k}\|_{L^q_{t,x}}
\end{align*}
By Minkowski inequality, the above quantity can be controlled by
\begin{align*}&\Big\{\sum\limits_{j=0}^{k-1}\bigg(\int_{I}\Big(\int_{\mathbb{R}^d}\frac{\|u_{(n),R}^{j}(x-y)-u_{(n),R}^{j}(x)
\|_{L^q}^2}{|y|}\frac{dy}{|y|^d}\Big)
^{\frac{q}{2}}\bigg)^{\frac2q}\Big\}^{\frac12}+\|u_{(n),R}^{<k}\|_{L^q_{t,x}}\\
&\lesssim\Big(\sum\limits_{j=0}^{k-1}\|u_{(n),R}^j\|_{L^q(I;B^\frac12_{q,2})}^2\Big)^\frac12.
\end{align*}
This implies that
$$\|u_{(n),R}^{<k}\|_{L^q(I;B^{\frac12}_{q,2})}\lesssim\Big(\sum\limits_{j=0}^{k-1}\|u_{(n),R}^j\|_{L^q(I;B^\frac12_{q,2})}^2\Big)^\frac12.$$
Thus we obtain \eqref{equ3.9}.

Finally, we turn to prove \eqref{equ3.10}. After the smooth cut-off,
we have for large $n$
\begin{equation}\label{equati}\big|\sum\limits_{j=0}^{k-1}u_{(n),R}^j\big|^2=\sum\limits_{j=0}^{k-1}|u_{(n),R}^j|^2.\end{equation}
Hence by the triangle inequality and \eqref{equati}, we obtain for
large $n$
\begin{align*}
&\big\|(|x|^{-\gamma}*|u_{(n)}^{<k}|^2)u_{(n)}^{<k}-\sum\limits_{j<k}(|x|^{-\gamma}*|u_{(n)}^j|^2)u_{(n)}^j\big\|_{ST^*(I)}\\
\leq&\big\|(|x|^{-\gamma}*|u_{(n)}^{<k}|^2)u_{(n)}^{<k}-(|x|^{-\gamma}*|u_{(n),R}^{<k}|^2)u_{(n),R}^{<k}\big\|_{ST^*(I)}\\
&+\sum\limits_{j=0}^{k-1}\big\|(|x|^{-\gamma}*|u_{(n),R}^{j}|^2)u_{(n),R}^{j}-(|x|^{-\gamma}*|u_{(n)}^{j}|^2)u_{(n)}^{j}\big\|_{ST^*(I)}\\
&+\sum\limits_{0\leq j\neq
l<k}\big\|(|x|^{-\gamma}*|u_{(n),R}^j|^2)u_{(n),R}^l\big\|_{ST^*(I)},
\end{align*}
therefore it suffices to prove
\begin{equation}\label{equation1}\big\|(|x|^{-\gamma}*|u_{(n)}^{<k}|^2)u_{(n)}^{<k}-(|x|^{-\gamma}*|u_{(n),R}^{<k}|^2)u_{(n),R}^{<k}\big\|_{ST^*(I)}\rightarrow
0,\end{equation} and
\begin{equation}\label{equation2}\big\|(|x|^{-\gamma}*|u_{(n),R}^j|^2)u_{(n),R}^l\big\|_{ST^*(I)}\rightarrow
0,\ j\not=l.\end{equation} For \eqref{equation1}, using the triangle
inequality , we have
\begin{align*}
\text{LHS}\ \text{of}\
\eqref{equation1}\leq&\Big\|(|x|^{-\gamma}*|u_{(n)}^{<k}|^2)(u_{(n)}^{<k}-u_{(n),R}^{<k})\Big\|_{ST^*(I)}\\
&+\Big\|\big[|x|^{-\gamma}*
\big(|u_{(n)}^{<k}|^2-|u_{(n),R}^{<k}|^2\big)u_{(n),R}^{<k}\Big\|_{ST^*(I)}\\
\triangleq& I_1+I_2.
\end{align*}
For $I_1$, by \eqref{eq1.8}, we get
\begin{align}\nonumber
I_1\lesssim&\big\| v
\big\|_{[K](I)}\|u_{(n)}^{<k}\|_{[K](I)}^{\frac4d}\|u_{(n)}^{<k}\|_{L_t^\infty
L_x^2}^{\frac{2(d-2)}{d}}+\big\| u_{(n)}^{<k}
\big\|_{[K](I)}^{1+\frac2d}\|u_{(n)}^{<k}\|_{L_t^\infty
L_x^2}^{\frac{d-2}{d}}\|v\|_{[K](I)}^{\frac2d}\|v\|_{L_t^\infty
L_x^2}^{\frac{d-2}{d}}\\\nonumber &+\big\| v
\big\|_{[W](I)}\|u_{(n)}^{<k}\|_{L_t^\infty
\dot{H}^{1}_x}^{\frac{2(d-3)}{d-1}}\|u_{(n)}^{<k}\|_{[W](I)}
^{\frac{4}{d-1}}+\big\| u_{(n)}^{<k}
\big\|_{[W](I)}^{1+\frac{2}{d-1}}\|u_{(n)}^{<k}\|_{L_t^\infty
\dot{H}_x^1}^{\frac{d-3}{d-1}}\|v\|_{L_t^\infty
\dot{H}^{1}_x}^{\frac{d-3}{d-1}}\|v\|_{[W](I)}
^{\frac{2}{d-1}}\\\label{equation4} \rightarrow&
0,\quad\text{as}\quad R\rightarrow\infty,
\end{align}
where $v=u_{(n)}^{<k}-u_{(n),R}^{<k}$. Similarly, $I_2\rightarrow
0\quad\text{as}\quad R\rightarrow\infty.$

Now we turn to prove \eqref{equation2}. By the triangle inequality,
we get
\begin{align}\nonumber
\text{LHS}\ \text{of}\
\eqref{equation2}&\leq\|(V_1*|u_{(n),R}^j|^2)u_{(n),R}^l\|_{[W]^*(I)}+\|(V_2*|u_{(n),R}^j|^2)u_{(n),R}^l\|_{[K]^*(I)}\\\label{equation5}&\triangleq
I_3+I_4,
\end{align}
where $V_1=|x|^{-\gamma}\chi_{|x|\leq R},
V_2=|x|^{-\gamma}\chi_{|x|\geq R}$.

For $I_3$, by the compact support of the function $u_{(n),R}^j$ and
\eqref{eq1.4}, we obtain
\begin{align}\nonumber
I_3=&\big\|(V_1*|u_{(n),R}^j(x)|^2)\chi_{|x|\lesssim
R}(x)u_{(n),R}^l(x-(x_n^l-x_n^j))\big\|_{[W]^*(I)}\\\nonumber
\lesssim& \big\| \chi_{|x|\lesssim R}(x)u_{(n),R}^l(x-(x_n^l-x_n^j))
\big\|_{[W](I)}\|u_{(n),R}^j(x)\|_{L_t^\infty
\dot{H}^{1}_x}^{\frac{2(d-3)}{d-1}}\|u_{(n),R}^j(x)\|_{[W](I)}
^{\frac{4}{d-1}}\\\nonumber&+\big\| u_{(n),R}^j(x)
\big\|_{[W](I)}^{1+\frac{2}{d-1}}\|u_{(n),R}^j(x)\|_{L_t^\infty
\dot{H}_x^1}^{\frac{d-3}{d-1}}\|\chi_{|x|\lesssim
R}(x)u_{(n),R}^l(x-(x_n^l-x_n^j))\|_{L_t^\infty
\dot{H}^{1}_x}^{\frac{d-3}{d-1}}\\\nonumber&\qquad\times\|\chi_{|x|\lesssim
R}(x)u_{(n),R}^l(x-(x_n^l-x_n^j))\|_{[W](I)}
^{\frac{2}{d-1}}\\\label{equation6} \rightarrow&\ 0\quad\text{as}\quad
n\rightarrow\infty,
\end{align}

For $I_4$, by the same argument as \eqref{eq1.3}, and
$\|V_2\|_{L^{\frac{d}{2}}}=R^{(1-\frac{\gamma}{2})d}$, we have
\begin{align}\nonumber
I_4=&\|(V_2*|u_{(n),R}^j|^2)u_{(n),R}^l\|_{[K]^*(I)}\\\nonumber
\lesssim& R^{(1-\frac{\gamma}{2})d}\Big\{\big\| u_{(n),R}^l
\big\|_{[K](I)}\|u\|_{[K](I)}^{\frac4d}\|u_{(n),R}^j\|_{L_t^\infty
L_x^2}^{\frac{2(d-2)}{d}}\\\nonumber&\qquad\qquad+\big\| u_{(n),R}^j
\big\|_{[K](I)}^{1+\frac2d}\|u_{(n),R}^j\|_{L_t^\infty
L_x^2}^{\frac{d-2}{d}}\|u_{(n),R}^l\|_{[K](I)}^{\frac2d}\|u_{(n),R}^l\|_{L_t^\infty
L_x^2}^{\frac{d-2}{d}}\Big\}\\\label{equation7}\rightarrow& 0\quad
\text{as}\quad R\rightarrow \infty.
\end{align}

This concludes the proof.

\end{proof}
After this preliminaries, we now show that
$\vec{u}_{(n)}^{<k}+\vec{w}_n^k$ is a good approximation for
$\vec{u}_n$ provided that each nonlinear profile has finite global
Strichartz norm.
\begin{lemma}\label{precldes}
Assume $\mu=1$. Let $u_n$ be a sequence of local solutions of
\eqref{equ1} around $t=0$ satisfying
$\varlimsup\limits_{n\rightarrow\infty}E(u_n,\dot{u}_n)<+\infty.$
Suppose that in its nonlinear profile decomposition
\eqref{nonlineard}, every nonlinear profile $\vec{U}^j_\infty$ has
finite global Strichartz and energy norms, i.e.
\begin{equation}\label{equ3.11}
\|\Re\langle\nabla\rangle^{-1}\vec{U}^j_\infty\|_{ST(\mathbb{R})}+\|\vec{U}^j_\infty\|_{L_t^\infty
L_x^2(\mathbb{R})}<\infty.
\end{equation}
Then $u_n$ is bounded for large $n$ in the Strichartz and the energy
norms, i.e.
\begin{equation}\label{equ3.12}
\varlimsup\limits_{n\rightarrow\infty}\|u_n\|_{ST(\mathbb{R})}+\|\vec{u}_n\|_{L_t^\infty
L_x^2(\mathbb{R})}<+\infty.
\end{equation}

Moreover, assume $\mu=-1$ and let $u_n$ be a sequence of local
solutions of \eqref{equ1} around $t=0$ in $\mathcal{A}^+$ satisfying
$\varlimsup\limits_{n\rightarrow\infty}E(u_n,\dot{u}_n)<E(\mathcal{W},0).$
Then the above results also hold true.
\end{lemma}

\begin{proof}
We only need to verify the condition of Lemma \ref{long}. For this purpose, we always use the fact that
$u_{(n)}^{<k}+w_n^k$ satisfies that
\begin{align*}
(\partial_{tt}+1-\Delta)\big(u_{(n)}^{<k}+w_n^k\big)=-&f(u_{(n)}^{<k}+w_n^k)+\big(f(u_{(n)}^{<k}+w_n^k)-f(u_{(n)}^{<k})\big)\\&+\big(f(u_{(n)}^{<k})-\sum
\limits_{j=0}^{k-1}f(u_{(n)}^j)\big).
\end{align*}

First, by the definition of the nonlinear concentrating wave
$u_{(n)}^j$ and Remark \ref{remark5.2}, we have
\begin{align}\nonumber
\big\|\big(\vec{u}_{(n)}^{<k}(0)+\vec{w}_n^k(0)\big)-\vec{u}_n(0)\big\|_{L^2_x}&\leq\sum\limits_{j=0}^{k-1}
\big\|\vec{u}_{(n)}^j(0)-\vec{u}_n^j(0)\big\|_{L^2_x}\rightarrow0,
\end{align}
as $n\rightarrow +\infty.$

Next, by the linear profile decomposition in Lemma \ref{lem3.1}, we
get
\begin{equation}\label{ee5}
\|\vec{u}_n(0)\|_{L^2}^2=\|\vec{v}_n(0)\|_{L^2}^2\geq\sum\limits_{j=0}^{k-1}\|\vec{v}_n^j(0)\|_{L^2}^2+o_n(1)
=\sum\limits_{j=0}^{k-1}\|\vec{u}_{(n)}^j(0)\|_{L^2}^2+o_n(1).
\end{equation}
Hence except for a finite set $J\subset\mathbb{N}$, the energy of
$u_{(n)}^j$ with $j\not\in J$ is smaller than the iteration
threshold (the small data  scattering in Lemma \ref{small}), and so
$$\|u_{(n)}^j\|_{ST(\mathbb{R})}\lesssim\|\vec{u}_{(n)}^j(0)\|_{L^2_x},\quad
j\not\in J.$$ This together
\eqref{equ3.8}, \eqref{equ3.9}, \eqref{equ3.11} and \eqref{ee5} yields that for any finite interval $I$
\begin{align}\nonumber
\sup\limits_{k}\varlimsup\limits_{n\rightarrow\infty}\|u_{(n)}^{<k}\|_{ST(I)}^2&\lesssim\sum\limits_{j\in
J}\|u_{(n)}^{j}\|_{ST(\mathbb{R})}^2+\sum\limits_{j\not\in
J}\|u_{(n)}^j\|_{ST(\mathbb{R})}^2\\\label{ee6} &\lesssim\sum\limits_{j\in
J}\|\Re\langle\nabla\rangle^{-1}\vec{U}_\infty^j\|_{ST(\mathbb{R})}^2+\varlimsup\limits_{n\rightarrow\infty}\|\vec{u}_n(0)\|_{L^2}^2<+\infty.
\end{align}
This together with the Strichartz estimate for $w_n^k$ implies that
$$\sup\limits_{k}\varlimsup\limits_{n\rightarrow\infty}\|u_{(n)}^{<k}+w_n^j\|_{ST(I)}<+\infty.$$
At last, by Lemma \ref{nonlineard} and Lemma \ref{lem3.3}, we have
$$\|f(u_{(n)}^{<k}+w_n^k)-f(u_{(n)}^{<k})\|_{ST^*(I)}\rightarrow
0,$$ and
$$\big\|f(u_{(n)}^{<k})-\sum
\limits_{j=0}^{k-1}f(u_{(n)}^j)\big\|_{ST^*(I)}\rightarrow 0,$$ as
$n\rightarrow+\infty.$ Therefore, by Lemma \ref{long}, we can obtain
the desired result.

\end{proof}

\section{Concentration Compactness}
\setcounter{section}{6}\setcounter{equation}{0} By the profile
decomposition in the previous section and the stability theory , we argue in this section
that if the scattering result does not hold, then there must exist a
minimal energy solution with some good compactness properties. This is the object of the
following proposition.
\begin{proposition}\label{prop4.1}
Let $\mu=1$. Suppose that  $E_{max}^+<+\infty$. Then there exists a
global solution $u_c$ of \eqref{equ1} satisfying
\begin{equation}\label{critical1}
E(u_c)=E_{max},\quad \|u_c\|_{ST(\mathbb{R})}=+\infty.
\end{equation}
Moreover, there exists $c(t):\mathbb{R}\rightarrow\mathbb{R}^d$,
such that $K=\{(u_c,\dot{u}_c)(t,x-c(t))\ \big|\ t\in\mathbb{R}^+\}$
is precompact in $H^1\times L^2$. Besides, one can assume that
$c(t)$ is $C^1$ and satisfies
\begin{equation}\label{smooth}
|\dot{c}(t)|\lesssim_{u_c} 1
\end{equation}
uniformly in $t$.

Furthermore, let $\mu=-1$, and suppose that $E_{max}^-<E(\mathcal{W},0)$. Then
the above results also hold true and
$(u_c,\dot{u}_c)\in\mathcal{A}^+.$
\end{proposition}

\begin{proof}
By the definition of $E_{max}$, we can choose a sequence
$\{u_n(t)\}$ such that \begin{equation}\label{defocusing}
E(u_n,\dot{u}_n)\rightarrow E_{max}^+,\ \text{and}\
\|u_n\|_{ST(\mathbb{R})}\rightarrow\infty,\ \text{as}\
n\rightarrow\infty,\ \text{(defocusing)}\end{equation} or
\begin{align}\label{focusing}
\begin{cases}
E(u_n,\dot{u}_n)\rightarrow E_{max}^-,\ \text{and}\
\|u_n\|_{ST(\mathbb{R})}\rightarrow\infty,\ \text{as}\
n\rightarrow\infty,\\
\|u_n\|_2^2+\|\nabla u_n\|_2^2<\|\mathcal{W}\|_2^2+\|\nabla \mathcal{W}\|_2^2,
\end{cases}\text{(focusing)}
\end{align}
then $u_n$ is global by Proposition \ref{pro23} in the defocusing
case and by Corollary \ref{global} in the focusing case. Now we
consider the linear and nonlinear profile decompositions of $u_n$,
using Lemma \ref{lem3.1},
\begin{align}\nonumber
e^{it\langle\nabla\rangle}\vec{u}_n(0)=\sum\limits_{j=0}^{k-1}
\vec{v}^j_n+\vec{\omega}_n^k,\
\vec{v}^j_n=e^{i\langle\nabla\rangle(t-t^j_n)}\varphi^j(x-x^j_n),\\
u_{(n)}^{<k}=\sum\limits_{j=0}^{k-1}u_{(n)}^j,\
\vec{u}^j_{(n)}(t,x)=\vec{U}^j_\infty(t-t_n,x-x_n),\\\nonumber
\|\vec{v}_n^j(0)-\vec{u}_{(n)}^j(0)\|_{L^2_x}\rightarrow 0, \ as\
n\rightarrow\infty.
\end{align}
Lemma \ref{precldes} precludes that all the nonlinear profiles
$\vec{U}_\infty^j$ have finite global Strichartz norm. On the other
hand, every solution of \eqref{equ1} with energy less than $E_{max}$
has global finite Strichartz norm by the definition of $E_{max}$.
Hence by \eqref{equ3.5}, we deduce that there is only one profile,
i.e. $K=1,$ and so for large $n$
\begin{equation}
\tilde{E}(\vec{u}_{(n)}^0)=E_{max},\
\|U_\infty^0\|_{ST(\mathbb{R})}=\infty,\
\lim\limits_{n\rightarrow\infty}\|\vec{\omega}_n^1\|_{L_t^\infty
L^2_x}=0.
\end{equation}
Hence we obtain \eqref{critical1} for $u_c=U_\infty^0$, and also
there exist a sequence $(t_n,x_n)\in \mathbb{R}\times\mathbb{R}^d$
and $\phi\in L^2(\mathbb{R}^d)$ such that along some subsequence,
\begin{equation}\label{sequa}
\|\vec{u}_n(0,x)-e^{-it_n\langle\nabla\rangle}\phi(x-x_n)\|_{L^2_x}\rightarrow
0, \ n\rightarrow\infty.\end{equation} Moreover, since \eqref{equ1}
is symmetric in $t$, we may assume that
\begin{equation}\label{forward}
\|u_c\|_{ST(0,+\infty)}=+\infty.
\end{equation}
Now we need only find $c(t)$ satisfying the right properties. The
proof of \cite{Pa2} can be adapted verbatim, but we give a sketch
for the sake of completeness. For $(u,v)\in H^1\times L^2$ and
$y\in\mathbb{R}^d$, we define
\begin{align*}
E_0(u,v,y,R)&=\int_{|x-y|\leq R}\big(|v(x)|^2+|\nabla
u(x)|^2+|u(x)|^2\big)dx,\\
\lambda(u,v,R)&=\sup\limits_{y\in\mathbb{R}^d}E_0(u,v,y,R),\\
\rho(u,v,\delta)&=\inf\big\{R:\
\lambda(u,v,R)>(1-\delta)E_0(u,v)\big\}.
\end{align*}
We claim that for all fixed $\delta>0$,
$\rho(u_c(t),\dot{u}_c(t),\delta)$ remains bounded. In fact, if this
were not true, there would exist a sequence of times $\{t_n\}_n$
such that we have for all $n$ and all $y$, that
$E_0(u_c(t_n),\dot{u}_c(t_n),y,n)\leq(1-\delta)E_0(u_c(t_n),\dot{u}_c(t_n)).$
But the sequence
$\big\{(u_n(t),\dot{u}_n(t))=(u_c(t+t_n),\dot{u}_c(t+t_n))\big\}$
satisfies the hypothesis \eqref{defocusing} (or \eqref{focusing}),
hence by \eqref{sequa}, there exists a sequence $(t_n',Y_n)\in
\mathbb{R}\times\mathbb{R}^d$ and $\phi\in L^2(\mathbb{R}^d)$ such
that, up to a subsequence,
\begin{equation}\label{seq}
\|\vec{u}_n(0,x)-e^{-it_n'\langle\nabla\rangle}\phi(x-Y_n)\|_{L^2_x}\rightarrow
0, \ as\ n\rightarrow\infty.\end{equation} And so
\begin{equation}\label{equ7.2}
\|\vec{u}_c(t_n,x)-e^{-it_n'\langle\nabla\rangle}\phi(x-Y_n)\|_{L^2_x}\rightarrow
0, \ as\ n\rightarrow\infty.\end{equation} Now we claim that
$\{t_n'\}$ is bounded. Indeed, if $t_n^\prime\rightarrow-\infty$,
then by triangle inequality and Strichartz estimate, we have
\begin{align*}
\|\langle\nabla\rangle^{-1}e^{it\langle\nabla\rangle}\vec{u}_c(t_n)\|_{ST(0,\infty)}\leq&
\|\langle\nabla\rangle^{-1}e^{it\langle\nabla\rangle}\big(\vec{u}_c(t_n)-e^{-it_n^\prime\langle\nabla\rangle}\phi(x-x_n^\prime)\big)\|_{ST(0,\infty)}\\
&\ +\|\langle\nabla\rangle^{-1}e^{i(t-t_n^\prime)\langle\nabla\rangle}\phi(x-x_n^\prime)\|_{ST(0,\infty)}\\
\lesssim&\|\vec{u}_c(t_n,x)-e^{-it_n^\prime\langle\nabla\rangle}\phi(x-x_n^\prime)\|_{L^2_x}+\|\langle\nabla\rangle^{-1}e^{it\langle\nabla\rangle}\phi\|_{ST(-t_n^\prime,\infty)}\\
\rightarrow&\ 0,\ \text{as}\ n\rightarrow \infty,
\end{align*}
so that we can solve \eqref{equ1} of $u_c$ for $t>t_n$ with large
$n$ globally by iteration with small Strichartz norms, contradicting
with \eqref{forward}.

If $t_n^\prime\rightarrow+\infty$, by the similar argument, we get
$$\|\langle\nabla\rangle^{-1}
e^{it\langle\nabla\rangle}\vec{u}_c(t_n)\|_{ST(-\infty,0)}=\|\langle\nabla\rangle^{-1}e^{it\langle\nabla\rangle}\phi\|_{ST(-\infty,-t_n^\prime)}+o(1)\rightarrow
0,\ \text{as}\ n\rightarrow \infty.$$ And so by Lemma \ref{small},
we can solve \eqref{equ2} of $u$ for $t<t_n$ for large $n$ with
$$\|u_c(t)\|_{ST(-\infty,t_n)}\lesssim\|\langle\nabla\rangle^{-1}
e^{i(t-t_n)\langle\nabla\rangle}\vec{u}_c(t_n)\|_{ST(-\infty,t_n)}=\|\langle\nabla\rangle^{-1}
e^{it\langle\nabla\rangle}\vec{u}_c(t_n)\|_{ST(-\infty,0)}.$$ This
implies $u_c=0$ by taking the limit, which contradicts with
\eqref{forward}. Therefore $t_n^\prime$ is bounded, which means that $\{t_n'\}$ is
precompact.

And so  by \eqref{equ7.2}, there exists $(w_0,w_1)\in
H^1\times L^2$, up to a subsequence,
\begin{equation}\label{cauchy1}
\|(u_c(t_n,x+Y_n),\dot{u}_c(t_n,x+Y_n))-(w_0,w_1)\|_{H^1\times
L^2}\rightarrow 0, \ as\ n\rightarrow\infty.
\end{equation}
Consequently,
$$E_0(w_0,w_1,y,R)=\lim\limits_{n\rightarrow+\infty}E_0(u_c(t_n,x+Y_n),\dot{u}_c(t_n,x+Y_n),y,R)\leq(1-\delta)E_0(w_0,w_1)$$
for all $R$. Thus $E_0(w_0,w_1)\leq(1-\delta)E_0(w_0,w_1).$ This
gives $E_0(w_0,w_1)=0$, which contradicts the fact that
$E(w_0,w_1)=\lim\limits_{n\rightarrow+\infty}E(u_c(t_n),\dot{u}_c(t_n))=E_{max}.$
Consequently, there exists an decreasing function $R$ such that
$\rho(u_c(t),\dot{u}_c(t),\delta)<R(\delta)$ for all $t\geq0.$

A similar proof shows that there exists $\kappa(\delta)>0$ such that
for all $t\geq0$,
\begin{equation}\label{kappa}
\lambda(u_c(t),\dot{u}_c(t),R(\delta))>\kappa(\delta).
\end{equation}
We choose $\delta$ to be small such that $\delta<\frac{1}{24}$ and
$\sqrt{\delta}<\frac{\kappa(\delta)}{8E_{max}}$, and let $c(t)$ be such
that
$$\lambda(u_c(t),\dot{u}_c(t),R(\delta))=E_0(u_c(t),\dot{u}_c(t),-c(t),R(\delta)).$$
We claim that the set $K=\big\{(u_c,\dot{u}_c)(t,x-c(t))\ \big|\
t\in\mathbb{R}^+\big\}$ is precompact in $H^1\times L^2$. Suppose it
were not true, then there would exist $\varepsilon>0$ and a sequence
of times $t_i$ such that
\begin{equation}\label{cauchy}
E_0\big(u_c(t_i,x-c(t_i))-u_c(t_j,x-c(t_j)),\dot{u}_c(t_i,x-c(t_i))-\dot{u}_c(t_j,x-c(t_j))\big)>\varepsilon
\end{equation}
for all $i\neq j.$ Using the same argument as \eqref{cauchy1}, we
obtain that there exists a sequence $\{Y_k\}_k$ and $(w_0,w_1)\in
H^1\times L^2$ such that, up to a subsequence,
$$(U(t_i),U_t(t_i))=(u_c(t_i,x-(y(t_i)-Y_i)),\dot{u}_c(t_i,x-(y(t_i)-Y_i)))\rightarrow
(w_0,w_1)\ in\ H^1\times L^2, $$ as $i\rightarrow\infty.$ In
particular, $(U(t_i),U_t(t_i))$ is a Cauchy sequence. Let $i_0$ be
such that for all $j\geq i_0$, there holds that
\begin{equation}\label{cont}
E_0(U(t_i)-U(t_j),U_t(t_i)-U_t(t_j))<\frac{\kappa(\delta)}{4},
\end{equation}
and suppose that there exists a subsequence such that
$|Y_k-Y_{i_0}|\rightarrow+\infty$ as $k\rightarrow+\infty$. Then,
for $|Y_k-Y_{i_0}|>2R(\delta)$, we have
\begin{align*}
&E_0(U(t_j)-U(t_{i_0}),U_t(t_j)-U_t(t_{i_0}))\\
=&E_0(U(t_j),U_t(t_j))+E_0(U(t_{i_0}),U_t(t_{i_0}))\\
&-2\int_{\mathbb{R}^d}\big(U(t_{i_0})U_t(t_j)+\nabla
U(t_j)\nabla U(t_{i_0})+U(t_j)U(t_{i_0})\big)\\
\geq&2\kappa(\delta)-2\int_{|x-Y_j|\leq
R(\delta)}\big(U(t_{i_0})U_t(t_j)+\nabla
U(t_j)\nabla U(t_{i_0})+U(t_j)U(t_{i_0})\big)\\
&-2\int_{|x-Y_j|\geq R(\delta)}\big(U(t_{i_0})U_t(t_j)+\nabla
U(t_j)\nabla U(t_{i_0})+U(t_j)U(t_{i_0})\big)\\
\geq&2\kappa(\delta)-2\sqrt{E_0(U(t_j),U_t(t_j))}\sqrt{\delta
E_0(U(t_{i_0}),U_t(t_{i_0}))}\\
&-2\sqrt{E_0(U(t_{i_0}),U_t(t_{i_0}))}\sqrt{\delta
E_0(U(t_j),U_t(t_j))}\\
\geq&2\kappa(\delta)-8\sqrt{\delta}E_{max}\geq\kappa(\delta),
\end{align*}
this contradicts with \eqref{cont}, where we use the fact that $E_0(u_c,\dot{u}_c)\leq2E(u_c,\dot{u}_c)$ in the third inequality. Thus, the sequence
$\{Y_k\}_k$ remains bounded. Therefore, up to a subsequence, we can
assume that $Y_k\rightarrow Y_*$. This fact implies that
$$\big(U(t_n,x+Y_*),U_t(t_n,x+Y_*)\big)=\big(u_c(t_n,\cdot-y(t_n)-(Y_n-Y_*)),\dot{u}_c(t_n,\cdot-y(t_n)-(Y_n-Y_*))\big)$$
is a Cauchy sequence, which contradicts \eqref{cauchy}. And so the
set $K=\big\{(u_c,\dot{u}_c)(t,x-c(t))\ \big|\
t\in\mathbb{R}^+\big\}$ is precompact in $H^1\times L^2$.

It only remains to prove \eqref{smooth}. By the precompactness of
$K$, and the continuity of the flow, there exists $s_0>0$ such that
for every solution $u$ of \eqref{equ1} with initial data
$(w_0,w_1)\in K$, there holds that
\begin{align*}
E_0(u_c(s),\dot{u}_c(s),0,2R(\delta))&\geq(1-\delta)E_0(u_c(0),\dot{u}_c(0),0,R(\delta)),\\
E_0(u_c(s),\dot{u}_c(s))&\leq(1-\delta)^{-1}E_0(u_c(0),\dot{u}_c(0))
\end{align*}
for every time $s$ such that $|s|\leq s_0$. In particular,
$$E_0(u_c(s),\dot{u}_c(s),0,2R(\delta))\geq(1-\delta)^3E_0(u_c(s),\dot{u}_c(s)).$$
This implies that
$E_0(u_c(s),\dot{u}_c(s),Y,R(\delta))<(1-\delta)E_0(u_c(s),\dot{u}_c(s))$
when $|Y|>3R(\delta).$ Consequently, for all $t\geq0$, for all
sufficiently small $s\leq s_0$, there holds that $|c(t)-c(t+s)|\leq
6R(\delta)$. Now, let $t_j=js_0$ for $j\in\mathbb{N}$, and let
$\tilde{c}(t)$ be a smooth function such that
$\tilde{c}(t_j)=c(t_j)$, and
$|\tilde{c}'(t)|\leq8R(\delta)s_0^{-1}.$ Then
$|c(t)-\tilde{c}(t)|\leq14R(\delta)\lesssim_{u_c}1,$ hence
$\big\{(u_c(t,x-\tilde{c}(t),\dot{u}_c(t,x-\tilde{c}(t)\big|
t\in\mathbb{R}^+\big\}$ also is precompact in $H^1\times L^2$ and
replacing $c(t)$ by $\tilde{c}(t)$, we obtain \eqref{smooth}. This
concludes the proof.
\end{proof}

As a consequence of the above proposition and
Hardy-Littlewood-Sobolev inequality, we have
\begin{corollary}\label{quick}
If we denote$$E_{R,c}=\int_{|x-c|\geq R}\big(|u|^2+|\nabla
u|^2+|\dot{u}|^2\big)dx+\iint\limits_{\substack{|x-c|\geq R \\ y\in
\mathbb{R}^d}}\frac{|u(t,x)|^2|u(t,y)|^2}{|x-y|^\gamma}dx dy,$$ then
for any $\eta>0$, there exists $R(\eta)>0$ such that
$$E_{R(\eta),c(t)}\leq\eta E(u,\dot{u}),\ for\ any\ t>0.$$
\end{corollary}
\begin{remark}\label{radial}
In particular, for the radial data $(u_0,u_1)\in H^1\times L^2$, by the same argument
as in \cite{KM,KV,MXZ}, we have $c(t)\equiv 0$, i.e.
$$E_{R(\eta),0}\leq\eta E(u,\dot{u}),\ for\ any\ t>0.$$
\end{remark}

The next corollary is the conclusion of this section.
\begin{corollary}\label{cor4.2}
Let $u$ be a nonlinear strong solution of \eqref{equ1} such that the
set $K$ defined in Proposition \ref{prop4.1} is precompact in
$H^1\times L^2$, and $E(u,\dot{u})\neq 0.$ Then there exists a
constant $\beta=\beta(\tau)>0$ such that, for all time $t,$ there
holds that
\begin{align}\label{big}
\int_t^{t+\tau}\iint_{\mathbb{R}^d\times\mathbb{R}^d}\frac{|x_2-y_2|^2}{|x-y|^{\gamma+2}}|u(s,x)|^2|u(s,y)|^2dx
dyds \geq\beta,\\
\int_t^{t+\tau}\iint_{\mathbb{R}^d\times\mathbb{R}^d}\frac{1}{|x-y|^{\gamma}}|u(s,x)|^2|u(s,y)|^2dx
dyds \geq\beta.
\end{align}
In particular, there holds that
$\int_0^{t}\iint_{\mathbb{R}^d\times\mathbb{R}^d}\frac{|x_2-y_2|^2}{|x-y|^{\gamma+2}}|u(t,x)|^2|u(t,y)|^2dx
dyds\gtrsim t.$
\end{corollary}

\begin{proof}
If \eqref{big} were not true, then there exists $\tau>0$ and a
sequence $\{t_k\}_k$ such that
\begin{equation}\label{nobound}
\int_t^{t+\tau}\iint_{\mathbb{R}^d\times\mathbb{R}^d}\frac{|x_2-y_2|^2}{|x-y|^{\gamma+2}}|u(t_k,x)|^2|u(t_k,y)|^2dx
dyds<\frac1k.
\end{equation}
Using the precompactness of $K$, we can extract a subsequence and
assume that $(u(t_k,\cdot-c(t_k)),\dot{u}(t_k,\cdot-c(t_k)))\to
(U_0,U_1)$ in $H^1\times L^2$. Let $U$ be the solution of
\eqref{equ1} with initial data $(U_0,U_1)$, then
$E(U,\dot{U})=E(u,\dot{u})\neq0$. Meanwhile, by \eqref {nobound} we
have
\begin{equation*}
\int_0^\tau\iint_{\mathbb{R}^d\times\mathbb{R}^d}\frac{|x_2-y_2|^2}{|x-y|^{\gamma+2}}|U(t,x)|^2|U(t,y)|^2dx
dydt=0.
\end{equation*}
Consequently, we have $U(t)=0$ a.e. for all $t\in(0,\tau)$, which
contradicts with $E(U_0,U_1)=E(U,\dot{U})=E(u,\dot{u})\neq0$. Hence,
we finish the proof.
\end{proof}

\section{Extinction of the critical element}
\setcounter{section}{7}\setcounter{equation}{0} In this section, we
prove that the critical solution constructed in Section 6 does not
exist, thus ensuring that $E_{max}^+ = +\infty$, and
$E_{max}^-=E(\mathcal{W},0)$. This implies Theorem \ref{theorem} and Theorem
 \ref{theorem1.2}.

\subsection{The defocusing case: $\mu=1$.}

\begin{proposition}\label{infty} Assume $d\geq 3,$ $\mu=1$, and $2<\gamma<\min\{4,d\}$, then $E_{max}^+=+\infty.$
\end{proposition}
\begin{proof}We use a Virial-type estimate in a direction orthogonal to the
Momentum vector. Up to relabeling the coordinates, we might assume
that $\hbox{Mom}(u)$ is parallel to the first coordinate. Thus we
have
\begin{equation}\label{equ4.1}\int_{\mathbb{R}^d}u_t(t,x)\partial_j u(t,x)dx=0,\quad \forall
j\geq 2.\end{equation}
 Let $\phi_R(x)=\phi(x/R)$ where
$\phi(x)$ is a nonnegative smooth radial function such that
supp$\phi\subseteq B(0,2)$ and $\phi\equiv1$ in $B(0,1)$. We define
the Virial action
$$I(t)=\int_{\mathbb{R}^d}z_2\phi_R(z)\partial_2u(t,x)u_t(t,x)dx,$$
where $z=x-c(t)$ and $z_2$ denotes the second component of
$z\in\mathbb{R}^d$. Integrating by parts we
get by \eqref{equ1}
\begin{align*}\nonumber
\partial_tI(t)=&\int_{\mathbb{R}^d}\partial_t(z_2\phi_R(z))\partial_2u(t,x)u_t(t,x)dx+\frac12\int_{\mathbb{R}^d}z_2\phi_R(z)\partial_2(u_t(x,t))^2dx\\\nonumber
&+\int_{\mathbb{R}^d}z_2\phi_R(z)\partial_2u(t,x)\big(\Delta
u-u-(|x|^{-\gamma}*|u|^2)u\big)dx\\
=&\frac12\int_{\mathbb{R}^d}\big(-|u_t|^2+|u|^2+|\nabla
u|^2+(|x|^{-\gamma}*|u|^2)|u|^2\big)dx-\int_{\mathbb{R}^d}|\partial_2u|^2dx\\
&+\dot{z}_2\int_{\mathbb{R}^d}u_t\partial_2udx-\frac{\gamma}{2}\int_{\mathbb{R}^d}z_2\phi_R(z)|u|^2(\frac{x_2}{|x|^{\gamma+2}}*|u|^2))dx\\
 &+\int_{|z|\geq R}\mathcal{O}_1(u)dx,
\end{align*}
where
\begin{align*}
\mathcal{O}_1(u)=&\frac12\big[\frac{z_2}{R}\phi_{R}'-(1-\phi_R(x))\big]\big[-|u_t|^2+|u|^2+|\nabla
u|^2+(|x|^{-\gamma}*|u|^2)|u|^2\big]\\
&-(c'(t)\cdot\nabla\phi_R)\frac{z_2}{R}\partial_2 u
u_t-c'_2(t)(1-\phi_R(z))\partial_2u u_t-(\nabla\phi_R\cdot\nabla
u)z_2\partial_2u,
\end{align*}
is bounded by a constant multiple of $(|u|^2+|\nabla
u|^2+|\partial_{t}u|^2)$ and is supported on the set $|z|\geq R$.
Besides, we define the equirepartition of energy action
$$J(t)=\int_{\mathbb{R}^d}\phi_R(z)u(t,x)u_t(t,x)dx.$$
Then
\begin{equation}
\partial_t J(t)=\int_{\mathbb{R}^d}\Big(|u_t|^2-|u|^2-|\nabla
u|^2-(|x|^{-\gamma}*|u|^2)|u|^2\Big)dx+\int_{|z|\geq
R}\mathcal{O}_2(u)dx,
\end{equation}
where
$$\mathcal{O}_2(u)=(1-\phi_R(z))\big[|u_t|^2-|u|^2-|\nabla
u|^2-(|x|^{-\gamma}*|u|^2)|u|^2\big]+\big(c^\prime(t)\cdot\nabla\phi_R\big)\frac{uu_t}{R}-\frac{u}{R}\nabla\phi_R\cdot\nabla
u,$$ has the same properties as $\mathcal{O}_1(u)$.

Considering $A(t)=I(t)+\frac12 J(t)$, we get
\begin{equation}\label{equ5.7}
|A(t)|\lesssim R E(u,\dot{u}),\ \text{for all time}\ t,
\end{equation}
 and
\begin{align*}
\partial_tA(t)=&-\int_{\mathbb{R}^d}|\partial_2u|^2dx-\frac{\gamma}{2}\iint\limits_{\mathbb{R}^d\times\mathbb{R}^d}\phi_R(x-c(t))(x_2-c_2(t))\frac{x_2-y_2}{|x-y|^{\gamma+2}}|u(t,x)|^2|u(t,y)|^2dxdy
\\&-\int_{|z|\geq R}(\mathcal{O}_1(u)+\frac12\mathcal{O}_2(u))dx.
\end{align*}
And so by symmetrization, $\partial_tA(t)$ can be rewritten as
\begin{align}\nonumber
-\partial_tA(t)=&\int_{\mathbb{R}^d}|\partial_2u|^2dx+\frac{\gamma}{4}\iint_{\mathbb{R}^d\times\mathbb{R}^d}\phi_R(z)\frac{|x_2-y_2|^2}{|x-y|^{\gamma+2}}|u(t,x)|^2|u(t,y)|^2dxdy\\\label{const2}
&+\frac{\gamma}{4}I_2+\int_{|z|\geq
R}(\mathcal{O}_1(u)+\mathcal{O}_2(u))dx,
\end{align}
where
\begin{align*}
I_2=\int\limits_{\mathbb{R}^d\times\mathbb{R}^d}&\big[(x_2-c_2(t))\phi_R(x-c(t))-(y_2-c_2(t))\phi_R(y-c(t))-(x_2-y_2)\big]\\
&\times\frac{x_2-y_2}{|x-y|^{\gamma+2}}|u(t,x)|^2|u(t,y)|^2dxdy.
\end{align*}

We will show that $I_2$ constitute only a small fraction of
$E(u,u_t)$. First, by Corollary \ref{cor4.2}, we know that if $R$ is
sufficient large depending on $u$ and $\eta$, then
$$E_{R,c(t)}(u,u_t)\leq \eta E(u,u_t).$$

Let $\chi$ denote a smooth cutoff to the region
$|x-c(t)|\geq\frac{R}{2}$ such that $\nabla\chi$ is bounded by
$R^{-1}$ and supported where $|x-c(t)|\sim R$. In the region where
$|x-c(t)|\thicksim|y-c(t)|$, we have
$$|x-c(t)|\thicksim|y-c(t)|\gtrsim R,$$ since otherwise $I_2$
vanish. Moreover, note that
$$|(x_2-c_2(t))\phi(x-c(t))-(y_2-c_2(t))\phi(y-c(t))|\lesssim|x-y|,$$
we use H\"{o}lder inequality to control the contribution to $I_2$
from this regime by
\begin{align*}
\iint_{\mathbb{R}^d\times \mathbb{R}^d}\frac{|\chi u(t,x)|^2|\chi
u(t,y)|^2}{|x-y|^\gamma}dxdy \lesssim\|\chi
u\|_2^{4-\gamma}\|\nabla(\chi u)\|_2^\gamma\lesssim \eta^2.
\end{align*}
In the region where $|x-c(t)|\ll|y-c(t)|,$ we use the fact that
$$|x-c(t)|\ll|y-c(t)|\sim|x-y|\ \text{and}\ |y-c(t)|\gtrsim R$$
to estimate the contribution from this regime by
$$\iint_{\mathbb{R}^d\times
\mathbb{R}^d}\frac{1}{|x-y|^\gamma}|\chi u(t,y)|^2|u(t,x)|^2
dxdy\lesssim \eta.$$ The last line follows from the same computation
as the first case. Finally, since the remaining region
$|y-c(t)|\ll|x-c(t)|$ can be estimated in the same way, we conclude
that
$$I_2\lesssim \eta.$$

Chosen $\eta$ sufficiently small depending on $u$ and $R$
sufficiently large depending on $u$ and $\eta$, we obtain
\begin{equation}\label{contrad}
-\partial_tA(t)\geq\frac{\gamma}{4}\iint_{\mathbb{R}^d\times\mathbb{R}^d}\phi_R(z)\frac{|x_2-y_2|^2}{|x-y|^{\gamma+2}}|u(t,x)|^2|u(t,y)|^2dxdy-\eta
E(u,u_t).
\end{equation}
If $E_{max}^+<\infty$, integrating \eqref{contrad} from 0 to $T>0$
and using Corollary \ref{cor4.2}, we get that there exists
$\alpha=\alpha(1,u)>0$ such that
$$\int_0^T\iint_{\mathbb{R}^d\times\mathbb{R}^d}\frac{|x_2-y_2|^2}{|x-y|^{\gamma+2}}|u(s,x)|^2|u(s,y)|^2dx
dy ds\geq \alpha T,$$ for all $T>1.$ Thus $-A(t)\gtrsim T$ for large
$T$, which contradicts with \eqref{equ5.7}. Hence we have
$E_{max}^+=+\infty$, this concludes the proof of Proposition
\ref{infty}.
\end{proof}

\subsection{The focusing case: $\mu=-1$}
\begin{proposition}\label{infty1} Assume $d\geq 3,$ $\mu=-1$, and $2<\gamma<\min\{4,d\}$, then $E_{max}^-=E(\mathcal{W},0).$
\end{proposition}

\begin{proof}
 Let $\phi_R(x)=\phi(x/R)$ where
$\phi(x)$ is a nonnegative smooth radial function such that
supp$\phi\subseteq B(0,2)$ and $\phi=1$ in $B(0,1)$. We define the
Virial action
$$I(t)=\int_{\mathbb{R}^d}\phi_R(x)x\cdot\nabla u(t,x)u_t(t,x)dx.$$
 Using \eqref{equ1}, and integrating
by parts we get
\begin{align}\nonumber
\partial_tI(t)=&\frac12\int_{\mathbb{R}^d}\phi_R(x)x\cdot\nabla(u_t(x,t))^2dx+\int_{\mathbb{R}^d}\phi_R(x)x\cdot\nabla u\big(\Delta
u-u+(|x|^{-\gamma}*|u|^2)u\big)dx\\
=&\frac{d}{2}\int_{\mathbb{R}^d}\big(-|u_t|^2+|u|^2+|\nabla
u|^2-(|x|^{-\gamma}*|u|^2)|u|^2\big)dx-\int_{\mathbb{R}^d}|\nabla u|^2dx\\
&+\frac{\gamma}{2}\int_{\mathbb{R}^d}\phi_R(x)|u|^2x\cdot(\frac{x}{|x|^{\gamma+2}}*|u|^2))dx+\int_{|x|\geq
R}\mathcal{O}_1(u)dx
\end{align}
where
\begin{align*}
\mathcal{O}_1(u)=&\big[\frac12\frac{x}{R}\cdot\nabla\phi_{R}-\frac{d}{2}(1-\phi_R(x))\big]\big[-|u_t|^2+|u|^2+|\nabla
u|^2-(|x|^{-\gamma}*|u|^2)|u|^2\big]\\
&+\frac12\sum\limits_{j=1}^d\frac{x_j}{R}\partial_j\phi_R|\partial_ju|^2-(\nabla\phi_R\cdot\nabla
u)\frac{x}{R}\cdot\nabla u,
\end{align*}
is bounded by a constant multiple of $(|u|^2+|\nabla
u|^2+|\partial_{t}u|^2)$ and is supported on the set $|x|\geq R$.
Besides, we define the equirepartition of energy action
$$J(t)=\int_{\mathbb{R}^d}\phi_R(x)u(t,x)u_t(t,x)dx.$$
Then
\begin{equation}
\partial_t J(t)=\int_{\mathbb{R}^d}\big(|u_t|^2-|u|^2-|\nabla
u|^2+(|x|^{-\gamma}*|u|^2)|u|^2\big)dx+\int_{|x|\geq
R}\mathcal{O}_2(u)dx,
\end{equation}
where
$$\mathcal{O}_2(u)=(1-\phi_R(x))\big[|u_t|^2-|u|^2-|\nabla
u|^2-(|x|^{-\gamma}*|u|^2)|u|^2\big]+\frac{u}{R}\nabla\phi_R\cdot\nabla
u,$$ has the same properties as $\mathcal{O}_1$. Considering
$A(t)=I(t)+\frac{d}{2} J(t)$, by the definition of $I(t)$ and
$J(t)$, we get that
\begin{equation}\label{equ5.9}
|A(t)|\lesssim R E(u,\dot{u}),
\end{equation}
for all time $t$ and
\begin{align*}
\partial_tA(t)=&-\int_{\mathbb{R}^d}|\nabla u|^2dx+\frac{\gamma}{2}\iint_{\mathbb{R}^d\times \mathbb{R}^d}\phi_R(x)\frac{x\cdot (x-y)}{|x-y|^{\gamma+2}}|u(t,x)|^2|u(t,y)|^2dxdy
\\&-\int_{|x|\geq R}(\mathcal{O}_1(u)+\mathcal{O}_2(u))dx,
\end{align*}
And so  by symmetrization, $\partial_tA(t)$ can be rewritten as
\begin{align*}
-\partial_tA(t)=&\int_{\mathbb{R}^d}|\nabla u|^2dx-\frac{\gamma}{4}\iint_{\mathbb{R}^d\times \mathbb{R}^d}\phi_R(x)\frac{|u(t,x)|^2|u(t,y)|^2}{|x-y|^{\gamma}}dxdy\\
&+\int_{|x|\geq R}(\mathcal{O}_1(u)+\mathcal{O}_2(u))dx\\
=&K_2(u)+\int_{|x|\geq R}(\mathcal{O}_1(u)+\mathcal{O}_2(u))dx,
\end{align*}
Using Proposition \ref{conservation} and Remark \ref{radial}, we
obtain
\begin{equation}\label{contrad1}
-\partial_tA(t)\geq\frac{\gamma}{4}\iint_{\mathbb{R}^d\times\mathbb{R}^d}\frac{|u(t,x)|^2|u(t,y)|^2}{|x-y|^{\gamma}}dxdy-\eta
E(u,u_t).
\end{equation}

 By Corollary
\ref{cor4.2}, we get that there exists $\alpha=\alpha(1,u)>0$ such
that
$$\int_0^T\iint_{\mathbb{R}^d\times\mathbb{R}^d}\frac{|u(s,x)|^2|u(s,y)|^2}{|x-y|^{\gamma}}dx
dy ds\geq \alpha T,$$ for all $T>1.$ Integrating the inequality
\eqref{contrad1} from 0 to $T>0$ and taking $\eta$ small, we obtain
$-A(t)\gtrsim T$ for large $T$, which contradicts with
\eqref{equ5.9}.

\end{proof}

\textbf{Acknowledgements} The authors are partly supported by the NSF
of China (No. 11171033). The
authors wish to thank Dr. Guixiang Xu, Haigen Wu and Junyong Zhang for stimulating discussion
about this problem.

\begin{center}

\end{center}

\end{document}